\documentclass{article}

\usepackage{microtype}
\usepackage[nottoc]{tocbibind}
\usepackage[utf8]{inputenc}
\usepackage{graphicx}
\usepackage{float}
\usepackage{bm}
\usepackage{amsthm}
\usepackage{amsmath}
\usepackage{amssymb}
\usepackage{mathrsfs}
\usepackage{dsfont}
\usepackage[all]{xy}
\usepackage{tikz-cd}
\usepackage{stmaryrd}
\usepackage{enumitem}
\usepackage{tabularx}
\usepackage{placeins}
\usepackage{pdflscape}
\usepackage{hyperref}

\theoremstyle{definition}
\newtheorem{Definition}{Definition}[subsection]

\theoremstyle{plain}
\newtheorem{Theorem}[Definition]{Theorem}

\theoremstyle{plain}

\newcounter{mainthm}

\newtheorem{maintheorem}[mainthm]{Theorem}

\theoremstyle{plain}
\newtheorem{Proposition}[Definition]{Proposition}

\theoremstyle{plain}
\newtheorem{Lemma}[Definition]{Lemma}

\theoremstyle{plain}
\newtheorem{Corollary}[Definition]{Corollary}

\theoremstyle{plain}
\newtheorem{Conjecture}[Definition]{Conjecture}

\theoremstyle{plain}

\theoremstyle{plain}

\theoremstyle{definition}
\newtheorem{Construction}[Definition]{Construction}

\theoremstyle{definition}
\newtheorem{Example}[Definition]{Example}

\theoremstyle{definition}

\theoremstyle{remark}
\newtheorem{Remark}[Definition]{Remark}

\title{On the Dualizability of Fusion 2-Categories}
\author{Thibault D. Décoppet}
\date{August 2024}

\begin{document}

\bibliographystyle{alpha}

\maketitle
    \hspace{1cm}
    \begin{abstract}
        Over an arbitrary field, we prove that the relative 2-Deligne tensor product of two separable module 2-categories over a compact semisimple tensor 2-category exists. This allows us to consider the Morita 4-category of compact semisimple tensor 2-categories, separable bimodule 2-categories, and their morphisms. Categorifying a result of \cite{DSPS13}, we prove that separable compact semisimple tensor 2-categories are fully dualizable objects therein. In particular, it then follows from the main theorem of \cite{D9} that, over an algebraically closed field of characteristic zero, every fusion 2-category is a fully dualizable object of the above Morita 4-category. We explain how this can be extended to any field of characteristic zero. Finally, we discuss the field theoretic interpretation of our results.
    \end{abstract}

\tableofcontents

\section*{Introduction}
\addcontentsline{toc}{section}{Introduction}

Over an algebraically closed field of characteristic zero, fusion 2-categories were introduced in \cite{DR} so as to define a state-sum invariant for oriented 4-manifolds. Categorifying the main result of \cite{DSPS13}, it was conjectured that every fusion 2-category is a fully dualizable object of an appropriate symmetric monoidal 4-category. In particular, thanks to the cobordism hypothesis \cite{BD,L}, this would readily imply that fusion 2-categories yield 4-dimensional framed fully extended topological field theories. An approach to the construction of the desired 4-category was given in \cite{JF} using the theory of higher condensations from \cite{GJF}. However, the latter theory has not yet been made completely rigorous:\ It depends on various assumptions regarding the behaviour of colimits in higher categories \cite[section 4.2]{GJF}.

We give a different approach to the construction of the sought-after 4-category. More precisely, we show that the relative 2-Deligne tensor product exists, thereby categorifying a result of \cite{ENO2} (see also \cite{DSPS14}). Our proof combines various properties of finite semisimple 2-categories and fusion 2-categories previously obtained in \cite{D1,D4,D7,D9}. With the relative 2-Deligne tensor product at our disposal, we may then appeal to the result of \cite{JFS} so as to construct a symmetric monoidal Morita 4-category $\mathbf{F2C}$, whose objects are fusion 2-categories, 1-morphisms are finite semisimple bimodule 2-categories, and higher morphisms are bimodule morphisms. Relatedly, the symmetric monoidal 4-category $\mathbf{BrFus}$ of braided fusion 1-categories, central fusion 1-categories, and their morphisms was studied in \cite{BJS}. We expect that the symmetric monoidal 4-categories $\mathbf{F2C}$ and $\mathbf{BrFus}$ are equivalent. Namely, it follows from \cite{D9} that these categories have the same equivalence classes of objects. The difficulty lies in exhibiting a symmetric monoidal 4-functor $\mathbf{BrFus}\rightarrow \mathbf{F2C}$. A sketch was given in \cite{JF}, but it does presuppose the currently unverified assumptions of \cite{GJF}.

Then, we prove that every fusion 2-category is a fully dualizable object of $\mathbf{F2C}$. Further, we show that every morphism has all possible adjoints. Our proof is inspired by that given in \cite{DSPS13} to show that every fusion 1-category is fully dualizable. In particular, our argument relies on the fact that every fusion 2-category over an algebraically closed field of characteristic zero is separable, which was established in \cite{D9}. A related dualizability result is obtained in \cite{BJS}, where it is shown that every braided fusion 1-category is a fully dualizable object of $\mathbf{BrFus}$. Indeed, provided that one was able to construct a symmetric monoidal 4-functor $\mathbf{BrFus}\rightarrow \mathbf{F2C}$, the results of \cite{BJS} and \cite{D9} would imply that fusion 2-categories over an algebraically closed field of characteristic zero are fully dualizable. In a slightly different direction, we wish to point out that we will subsequently work over an arbitrary perfect field, a level of generality at which the relation between these two approaches is, at present, unknown.

\subsection*{Results}

Working over an arbitrary field $\mathds{k}$, one can consider compact semisimple 2-categories. This notion was introduced in \cite{D5}, and categorifies that of a finite semisimple 1-category over $\mathds{k}$. Further, it specializes to that of a finite semisimple 2-category over algebraically closed fields of characteristic zero. We may then consider a compact semisimple tensor 2-category $\mathfrak{C}$, that is a monoidal compact semisimple 2-category whose objects have duals. This generalizes to an arbitrary field the concept of a multifusion 2-category. Further, there is a notion of compact semisimple $\mathfrak{C}$-module 2-category. We will mostly be interested in the following refinement of this notion:\ A compact semisimple module 2-category is said to be \textit{separable} if it is equivalent to the 2-category of modules over a separable algebra in $\mathfrak{C}$. These objects were studied in \cite{D8}, where an intrinsic characterization was established:\ Under mild technical assumptions, a compact semisimple $\mathfrak{C}$-module 2-category is separable if and only if the associated 2-category of $\mathfrak{C}$-module 2-endofunctors is compact semisimple. Our central result is the following theorem, which categorifies theorem 3.3 of \cite{DSPS14} (see also \cite[Section 3]{ENO2}).

\begin{maintheorem}[{Thm.\ \ref{thm:relative2Deligne}}]
Let $\mathfrak{C}$ be a compact semisimple 2-category. The relative 2-Deligne tensor product of any two separable compact semisimple $\mathfrak{C}$-module 2-categories exists.
\end{maintheorem}

\noindent If one works over an algebraically closed field of characteristic zero, it follows from \cite{D9} that the separability hypothesis is automatically satisfied.

Let us now assume that $\mathds{k}$ is a perfect field, this allows us to consider the symmetric monoidal 3-category $\mathbf{CSS2C}$ of compact semisimple 2-categories, linear 2-functors, 2-natural transformations, and modifications, whose monoidal structure is given by the 2-Deligne tensor product \cite{D3}. The symmetric monoidal 3-category $\mathbf{CSS2C}$ is equivalent to a sub-3-category of the symmetric monoidal Morita 3-category $\mathbf{TC}^{ss}$ of finite semisimple tensor 1-categories introduced in \cite{DSPS13}. Then, thanks to \cite{JFS} and the above theorem, we may consider the symmetric monoidal Morita 4-category $\mathbf{CSST2C}$ of compact semisimple tensor 2-categories, separable bimodule 2-categories, bimodule 2-functors, bimodule 2-natural transformations, and bimodule modifications. Its monoidal structure is given by the 2-Deligne tensor product $\boxtimes$ from \cite{D3}. Our main objective is to study the dualizability properties of the Morita 4-category $\mathbf{CSST2C}$. To this end, recall from \cite{DSPS13} that a finite semisimple tensor 1-category is fully dualizable in $\mathbf{TC}^{ss}$ if and only if it is separable, that is, has finite semisimple Drinfeld center. This implies readily that a compact semisimple 2-category is fully dualizable in $\mathbf{CSS2C}$ precisely if it is locally separable, that is, if its $End$-1-categories, which are finite semisimple tensor 1-categories, have finite semisimple Drinfeld centers. Then, by analogy with \cite{DSPS13}, it is natural to consider the following notion:\ A compact semisimple tensor 2-category $\mathfrak{C}$ is \textit{separable} if its underlying 2-category is locally separable, and the Drinfeld center of $\mathfrak{C}$ is a compact semisimple 2-category. Inspired by the techniques of \cite{DSPS13}, we derive the following result.

\begin{maintheorem}[Thm.\ \ref{thm:fulldualizability}]
Separable compact semisimple tensor 2-categories are fully dualizable objects of the symmetric monoidal 4-category $\mathbf{CSST2C}$.
\end{maintheorem}

\noindent In fact, we show that separability is automatically satisfied over any field of characteristic zero, extending a result of \cite{D9} over algebraically closed fields of characteristic zero.

\begin{maintheorem}[Cor.\ \ref{cor:separabilitycharzero}]
Over an arbitrary field of characteristic zero, every compact semisimple tensor 2-category is fully dualizable.
\end{maintheorem}

\noindent In particular, we bring a positive answer to the question raised in \cite{DR} of whether fusion 2-categories over an algebraically closed field of characteristic zero are fully dualizable. Finally, in section \ref{sub:TFTs} below, we discuss the field theoretic interpretation of these results.

\subsection*{Acknowledgments}

I would like to thank all the participants of the AIM SQuaRE ``Higher Symmetries of Fusion 2-Categories'' for encouragements and fruitful conversations. I am also grateful towards Christopher Douglas and Matthew Yu for comments on a draft of this manuscript. This work was supported in part by the Simons Collaboration on Global Categorical Symmetries.

\section{Preliminaries}

\subsection{2-Condensations}

The theory of 2-condensations as introduced in \cite{DR, GJF} is by now well-established. Nevertheless, as this notion will play a central part in the proof of our main theorem below, we do recall a number of key definitions and facts from \cite{GJF} and \cite{D1}.

\begin{Definition}
A 2-condensation in a 2-category $\mathfrak{C}$ consists of two objects $A$ and $B$, together with two 1-morphisms $f:A\rightarrow B$, $g:B\rightarrow A$, and two 2-morphisms $\phi:f\circ g\Rightarrow Id_B$, $\gamma:Id_B\Rightarrow f\circ g$, such that $\phi\cdot\gamma = Id_{Id_B}$.
\end{Definition}

The concept of a 2-condensation is a categorification of the notion of a split surjection. In particular, we may extract from the above definition a 2-categorical version of an idempotent, referred to as a 2-condensation monad.

\begin{Definition}
A 2-condensation monad in a 2-category $\mathfrak{C}$ consists of an object $A$, equipped with a 1-endomorphism $e:A\rightarrow A$, and two 2-morphisms $\mu:e\circ e\Rightarrow e$, $\delta:e\Rightarrow e\circ e$, such that $\mu$ is associative, $\delta$ is coassociative, the Frobenius relations hold, and $\mu\cdot\delta = Id_{e}$.
\end{Definition}

All the usual constructions involving idempotents have 2-categorical analogues. For instance, we say that a 2-condensation monad splits if it may be extended to a 2-condensation. If such an extension exists, it is essentially unique by theorem 2.3.2 of \cite{GJF}. Then, a 2-category $\mathfrak{C}$ is said to be Karoubi complete if it is locally idempotent complete, i.e.\ its $Hom$-1-categories have splittings for idempotents, and every 2-condensation monad splits. Furthermore, provided that the 2-category $\mathfrak{C}$ is locally idempotent complete, we can consider its Karoubi completion. That is, there exists a 2-functor $\iota:\mathfrak{C}\rightarrow Kar(\mathfrak{C})$ to a Karoubi complete 2-category \cite{GJF}. The 2-category $Kar(\mathfrak{C})$ has as objects 2-condensations monads in $\mathfrak{C}$, and as 1-morphisms 2-condensation bimodules (as introduced in \cite[section 2.3]{GJF}, see also \cite[Appendix A.2]{D1} for a spelled out definition). In fact, the 2-functor $\iota$ is 3-universal amongst 2-functors from $\mathfrak{C}$ to a Karoubi complete 2-category \cite{D1}. Let us remark that the assumption that $\mathfrak{C}$ be locally Karoubi complete is not restrictive:\ Every 2-category can be locally Karoubi completed. For later use, let us also record the following lemma, which follows from the observation that $\iota$ is fully faithful, that is, it induces equivalences on $Hom$-1-categories.

\begin{Lemma}\label{lem:splittings2condensationmonads}
Let $(A_i, e_i,\mu_i,\delta_i)$, $i=1,2$, be two 2-condensations monads in $\mathfrak{C}$ with splittings $(A_i, B_i, f_i, g_i,\phi_i,\gamma_i)$.
Then, $B_1$ and $B_2$ are equivalent objects in $\mathfrak{C}$ if and only if the two 2-condensation monads $(A_i, e_i,\mu_i,\delta_i)$ are equivalent as objects of $Kar(\mathfrak{C})$.
\end{Lemma}

Let us now assume that $\mathfrak{C}$ is an $R$-linear 2-category for some ring $R$ that is locally Cauchy complete, i.e.\ its $Hom$-1-categories have (finite) direct sums, and splittings for idempotents. Then it makes sense to talk about direct sums for the objects of $\mathfrak{C}$. We say that a linear 2-category is Cauchy complete is if it is locally Cauchy complete, it is Karoubi complete, and has direct sums for objects. There is also a notion of Cauchy completion $Cau(\mathfrak{D})$ for any locally Cauchy complete $R$-linear 2-category $\mathfrak{D}$. We refer the reader to \cite{GJF} for details (see also \cite{D1}).

\subsection{Compact Semisimple Tensor 2-Categories}

Over an algebraically closed field of characteristic zero, the notions of finite semisimple 2-category and fusion 2-category were introduced in \cite{DR}. As discussed in \cite{D5}, these notions can be generalized to an arbitrary field $\mathds{k}$ so as to give those of compact semisimple 2-category and compact semisimple tensor 2-category, which we now briefly review.

Let $\mathfrak{C}$ be a locally finite semisimple ($\mathds{k}$-linear) 2-category. An object $C$ of $\mathfrak{C}$ is called simple if $Id_C$ is a simple object of the finite semisimple 1-category $End_{\mathfrak{C}}(C)$. Two simple objects are called connected if there exists a non-zero 1-morphism between them. We write $\pi_0(\mathfrak{C})$ for the set of equivalences classes of simple objects for this relation; This is the set of connected components of $\mathfrak{C}$.

\begin{Definition}
A compact semisimple 2-category is a locally finite semisimple 2-category that is Cauchy complete, has right and left adjoints for 1-morphisms, and has finitely many connected components.
\end{Definition}

\noindent A finite semisimple 2-category is a compact semisimple 2-category that has finitely many equivalence classes of simple objects. Over an algebraically closed field of characteristic zero, it follows from \cite{DR} that every compact semisimple 2-category is finite. More generally, it was shown in \cite{D5} that, over an algebraically closed field of arbitrary characteristic or a real closed field, every compact semisimple 2-category is finite. On the other hand, over a perfect field that is neither algebraically closed nor real closed, every compact semisimple 2-category has infinitely many equivalence classes of simple objects.

Given a finite semisimple tensor 1-category $\mathcal{C}$, we can consider the 2-category $\mathbf{Mod}(\mathcal{C})$ of separable right $\mathcal{C}$-module 1-categories in the sense of \cite{DSPS13}. This is a compact semisimple 2-category, and every compact semisimple 2-category is of this form \cite{D5}. This generalizes a result of \cite{DR} over algebraically closed fields of characteristic zero.

Let us now assume that $\mathds{k}$ is a perfect field. Recall from \cite{DSPS13} that a finite semisimple tensor 1-category $\mathcal{C}$ is separable if its Drinfeld center $\mathcal{Z}(\mathcal{C})$ is finite semisimple. As is discussed extensively in this last reference, this condition is intimately related to fully dualizability. It is therefore not surprising that the following technical condition will appear repeatedly.

\begin{Definition}\label{def:locallyseparable}
A compact semisimple 2-category $\mathfrak{C}$ is locally separable if for every simple object $C$ in $\mathfrak{C}$, the finite semisimple tensor 1-category $End_{\mathfrak{C}}(C)$ is separable, that is, its Drinfeld center is finite semisimple.
\end{Definition}

\noindent In fact, it is enough to check this condition for one object in every connected component \cite{D5}.

Let $\mathds{k}$ be an arbitrary field. The following notion will be our main object of study.

\begin{Definition}
A compact semisimple tensor 2-category is a compact semisimple 2-category equipped with a monoidal structure such that every object admits a right and a left dual.
\end{Definition}

\noindent To every compact semisimple tensor 2-category $\mathfrak{C}$, there is an associated finite semisimple braided tensor 1-category $\Omega\mathfrak{C}:=End_{\mathfrak{C}}(I)$, the 1-category of 1-endo\-morphisms of the monoidal unit $I$ of $\mathfrak{C}$. If $\mathds{k}$ is an algebraically closed field or a real closed field, we will also use the term multifusion 2-category to refer to a compact semisimple tensor 2-category. A fusion 2-category is a multifusion 2-category whose monoidal unit $I$ is a simple object. We now give a few examples that will be relevant for our purposes. Many more may be found in \cite{DR} and \cite{DY23}.

\begin{Example}
If $\mathcal{B}$ is a finite semisimple braided tensor 1-category, then $\mathbf{Mod}(\mathcal{B})$ is a compact semisimple tensor 2-category with monoidal structure given by the relative Deligne tensor product $\boxtimes_{\mathcal{B}}$. For instance, we can take $\mathds{k}$ to be an algebraically closed field of characteristic different from $2$, and $\mathcal{B}=\mathbf{Vect}(\mathbb{Z}/4)$. The underlying 2-category of the fusion 2-category $\mathbf{Mod}(\mathbf{Vect}(\mathbb{Z}/4))$ is depicted below.

\begin{equation*}\begin{tikzcd}
& \mathbf{Vect}(\mathbb{Z}/2) \arrow["\mathbf{Vect}^{\kappa}(\mathbb{Z}/2\oplus \mathbb{Z}/2)"', loop, distance=2em, in=125, out=55] \arrow[ld, bend right] \arrow[rd, bend right = 10] &    \\
\mathbf{Vect}(\mathbb{Z}/4) \arrow["\mathbf{Vect}(\mathbb{Z}/4)"', loop, distance=4em, in=215, out=145] \arrow[rr, bend right=20] \arrow[ru, bend right =10] &  & \mathbf{Vect} \arrow["\mathbf{Vect}(\mathbb{Z}/4)"', loop, distance=3em, in=35, out=325] \arrow[lu, bend right] \arrow[ll, bend right=5]
\end{tikzcd}\end{equation*} We have used $\kappa$ to denote a specific 3-cocycle for $\mathbb{Z}/2\oplus\mathbb{Z}/2$ with coefficients in $\mathds{k}^{\times}$ (see \cite{Nai}). A sample of its fusion rules is given by \begin{equation*}\mathbf{Vect}(\mathbb{Z}/2)\boxtimes_{\mathbf{Vect}(\mathbb{Z}/4)}\mathbf{Vect}(\mathbb{Z}/2)=\boxplus_{i=1}^2\mathbf{Vect}(\mathbb{Z}/2),\end{equation*}
\begin{equation*}\mathbf{Vect}\boxtimes_{\mathbf{Vect}(\mathbb{Z}/4)}\mathbf{Vect}=\boxplus_{i=1}^4\mathbf{Vect}.\end{equation*}
\end{Example}

\begin{Example}
For simplicity, let us assume that $\mathds{k}$ is perfect. Given $G$ a finite group, and $\pi$ a 4-cocycle for $G$ with coefficients in $\mathds{k}^{\times}$, we can then consider the compact semisimple tensor 2-category $\mathbf{2Vect}^{\pi}(G)$ of $G$-graded finite semisimple 1-categories with pentagonator twisted by $\pi$. If $\mathds{k}$ is an algebraically closed field of characteristic zero, such fusion 2-categories were completely characterized in \cite{JFY} as those fusion 2-categories $\mathfrak{C}$ for which $\Omega\mathfrak{C}\simeq \mathbf{Vect}$. In fact, their argument does not depend on the characteristic, though it does require the base field to be algebraically closed. A bosonic strongly fusion 2-category is a fusion 2-category $\mathfrak{C}$ over an algebraically closed field (of characteristic zero) such that $\Omega\mathfrak{C}\simeq \mathbf{Vect}$. Likewise, a fermionic strongly fusion 2-category is a fusion 2-category $\mathfrak{C}$ such that $\Omega\mathfrak{C}\simeq \mathbf{SVect}$, the symmetric fusion 1-category of super vector spaces. The classification of fermionic strongly fusion 2-categories over an algebraically closed field of characteristic zero will be discussed in \cite{D11}. 
\end{Example}

\subsection{Separable Module 2-Categories}

Let $\mathfrak{C}$ be a monoidal 2-category with monoidal product $\Box$, and monoidal unit $I$. We can consider the notion of a (left) module 2-category over $\mathfrak{C}$. These were spelled out in \cite{D4}. Moreover, it was observed in proposition 2.2.8 of \cite{D4}, following a theorem of \cite{Gur}, that every pair consisting of a monoidal 2-category and a (left) module 2-category over it may be strictified. More precisely, we may assume without loss of generality that $\mathfrak{C}$ is a strict cubical monoidal 2-category, i.e.\ it satisfies the conditions of definition 2.26 of \cite{SP}, and the (left) module 2-category is as in the definition below.

\begin{Definition}\label{def:module2cat}
Let $\mathfrak{M}$ be a strict 2-category. A strict cubical left $\mathfrak{C}$-module 2-category structure on $\mathfrak{M}$ is a strict cubical 2-functor $\Box:\mathfrak{C}\times \mathfrak{M}\rightarrow \mathfrak{M}$ such that:
\begin{enumerate}
    \item The induced 2-functor $I\Box (-):\mathfrak{M}\rightarrow \mathfrak{M}$ is exactly the identity 2-functor,
    \item The two 2-functors \begin{equation*}\big((-)\Box (-)\big)\Box (-):\mathfrak{C}\times \mathfrak{C}\times\mathfrak{M}\rightarrow\mathfrak{M},\ \mathrm{and}\ (-)\Box \big((-)\Box (-)\big):\mathfrak{C}\times \mathfrak{C}\times\mathfrak{M}\rightarrow\mathfrak{M}\end{equation*} are equal on the nose.
\end{enumerate}
\end{Definition}

\noindent Let us also mention that this strictification procedure holds for any set of module 2-categories over a monoidal 2-category. In addition, there are notions of $\mathfrak{C}$-module 2-functors, $\mathfrak{C}$-module 2-natural transformations, and $\mathfrak{C}$-module modifications. There are also notions of right module 2-categories, as well as bimodule 2-categories, and maps between them. We refer the reader to \cite[section 2]{D4} for the precise definitions.

Let us now assume that $\mathds{k}$ is a perfect field, and that $\mathfrak{C}$ is a locally separable compact semisimple tensor 2-category. It was shown in \cite{D4} that every locally separable compact semisimple 2-category $\mathfrak{C}$-module 2-category is equivalent to the 2-category of modules over an algebra $A$ in $\mathfrak{C}$. However, just as in the decategorified setting \cite{DSPS13}, not all algebras in $\mathfrak{C}$ are well-behaved. More precisely, we are specifically interested in algebras for which the associated 2-categories of modules and bimodules are still compact semisimple. As explained in \cite{D7}, this is precisely achieved by the following refinement of the notion of an algebra.

\begin{Definition}\label{def:separable}
An algebra $A$ in a compact semisimple tensor 2-category is called separable if its multiplication 1-morphism $m:A\Box A\rightarrow A$ has a right adjoint $m^*$ as an $A$-$A$-bimodule 1-morphism, and the counit $\epsilon_A$ of this adjunction admits a splitting $\gamma_A$ as a 2-morphism of $A$-$A$-bimodules.\footnote{Every separable algebra is a 3-condensation monad in the sense of \cite{GJF}. In fact, we expect that every 3-condensation monad in a compact semisimple tensor 2-category is equivalent to a separable algebra in the 3-category of 3-condensation monads and condensation bimodules.}
\end{Definition}

\noindent We can now introduce a class of module 2-categories, which is particularly manageable.

\begin{Definition}
A compact semisimple left $\mathfrak{C}$-module 2-category $\mathfrak{M}$ is separable if it is equivalent to the left $\mathfrak{C}$-module 2-category $\mathbf{Mod}_{\mathfrak{C}}(A)$ of right $A$-modules over a separable algebra $A$ in $\mathfrak{C}$.
\end{Definition}

\noindent The above definition admits an obvious analogue for right $\mathfrak{C}$-module 2-categories. Further, we say that a $\mathfrak{C}$-$\mathfrak{D}$-bimodule 2-category is \textit{separable} if it is separable both as a left $\mathfrak{C}$-module 2-category and as a right $\mathfrak{D}$-module 2-category. Let us note that by proposition 5.3.1 of \cite{D8}, every separable compact semisimple module 2-category over a locally separable compact semisimple tensor 2-category is automatically locally separable.

For our current purposes, an important property of separable compact semisimple module 2-categories is theorem 5.3.2 of \cite{D8} recalled below.

\begin{Theorem}\label{thm:Moritadual}
Let $\mathfrak{M}$ be a separable compact semisimple left $\mathfrak{C}$-module 2-category. Then, $\mathbf{End}_{\mathfrak{C}}(\mathfrak{M})$, the monoidal 2-category of $\mathfrak{C}$-module 2-endofunctors on $\mathfrak{M}$, is a compact semisimple tensor 2-category.
\end{Theorem}

\noindent In fact, the converse of the above theorem also holds by proposition 5.1.6 of \cite{D8}. More precisely, a locally separable compact semisimple left $\mathfrak{C}$-module 2-category $\mathfrak{M}$ is separable if and only if $\mathbf{End}_{\mathfrak{C}}(\mathfrak{M})$ is compact semisimple.

\section{The Relative 2-Deligne Tensor Product}

\subsection{Balanced 2-Functors}

Let $\mathfrak{C}$ be a monoidal 2-category, $\mathfrak{M}$ a right $\mathfrak{C}$-module 2-category, and $\mathfrak{N}$ a left $\mathfrak{C}$-module 2-category. Without loss of generality, we assume that $\mathfrak{C}$, $\mathfrak{M}$, and $\mathfrak{N}$ are strict cubical in the sense of definition \ref{def:module2cat} above. We will use freely the variant of the graphical calculus of \cite{GS} given in \cite{D4}.

\begin{Definition}
Let $\mathfrak{A}$ be a strict 2-category. A $\mathfrak{C}$-balanced 2-functor $(\mathfrak{M},\mathfrak{N})\rightarrow \mathfrak{A}$ consists of:
\begin{enumerate}
    \item A 2-functor $\mathbf{B}:\mathfrak{M}\times\mathfrak{N}\rightarrow \mathfrak{A}$;
    \item An adjoint 2-natural equivalence $a^{\mathbf{B}}$ given on $M$ in $\mathfrak{M}$, $C$ in $\mathfrak{C}$, and $N$ in $\mathfrak{N}$ by \begin{equation*}a^{\mathbf{B}}:\mathbf{B}(M\Box C, N)\simeq \mathbf{B}(M,C\Box N)\end{equation*} with adjoint pseudo-inverse denoted by $(a^{\mathbf{B}})^{\bullet}$;
    \item Two invertible modifications $\omega^{\mathbf{B}}$ and $\gamma^{\mathbf{B}}$ given on $M$ in $\mathfrak{M}$, $C, D$ in $\mathfrak{C}$, and $N$ in $\mathfrak{N}$ by
\end{enumerate}

\begin{equation*}\begin{tikzcd}[sep=small]
                                                & {\mathbf{B}(MC,DN)} \arrow[rd, "a^{\mathbf{B}}"] \arrow[d, "\omega^{\mathbf{B}}", shorten <= 0.5ex,shorten >= 0.5ex, near start, Rightarrow] &            \\
{\mathbf{B}(MCD,N)} \arrow[ru, "a^{\mathbf{B}}"] \arrow[rr, "a^{\mathbf{B}}"'] & {}                                                 & {\mathbf{B}(M,CDN),}
\end{tikzcd}\end{equation*}

\begin{equation*}\begin{tikzcd}[sep=small]
{\mathbf{B}(MI,N)} \arrow[rr, "a^{\mathbf{B}}"] \arrow[rd,equal] & {} \arrow[d, "\gamma^{\mathbf{B}}", Rightarrow]                 & {\mathbf{B}(M,IN);} \\
                                                           & {\mathbf{B}(M,N)} \arrow[ru,equal] &                   
\end{tikzcd}\end{equation*}

satisfying:

\begin{enumerate}
    \item[a.] For every $M$ in $\mathfrak{M}$, $C,D,E$ in $\mathfrak{C}$, and $N$ in $\mathfrak{N}$, the equality
\end{enumerate}

\newlength{\prelim}
\settoheight{\prelim}{\includegraphics[width=45mm]{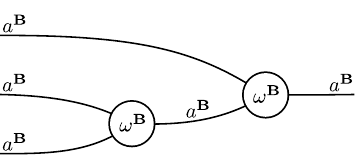}}

\begin{center}
\begin{tabular}{@{}ccc@{}}

\includegraphics[width=45mm]{pictures/balanced/balanced2pent1.pdf} & \raisebox{0.45\prelim}{$=$} &

\includegraphics[width=45mm]{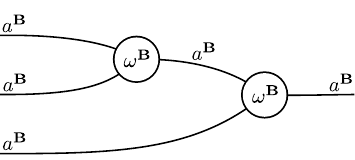}

\end{tabular}
\end{center}

\begin{enumerate}
    \item[] holds in $Hom_{\mathfrak{A}}(\mathbf{B}(M\Box C\Box D\Box E, N), \mathbf{B}(M,  C\Box D\Box E\Box N))$,

    \item[b.] For every $M$ in $\mathfrak{M}$, $C$ in $\mathfrak{C}$, and $N$ in $\mathfrak{N}$, the equality
    
\settoheight{\prelim}{\includegraphics[width=22.5mm]{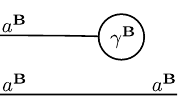}}

\begin{center}
\begin{tabular}{@{}ccc@{}}

\includegraphics[width=22.5mm]{pictures/balanced/balancedunitright1.pdf} & \raisebox{0.45\prelim}{$=$} &

\includegraphics[width=22.5mm]{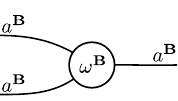}

\end{tabular}
\end{center}
    
    holds in $Hom_{\mathfrak{A}}(\mathbf{B}(M\Box C\Box I, N), \mathbf{B}(M, C\Box N))$,
    
    \item[c.] For every $M$ in $\mathfrak{M}$, $C$ in $\mathfrak{C}$, and $N$ in $\mathfrak{N}$, the equality

\settoheight{\prelim}{\includegraphics[width=22.5mm]{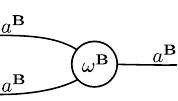}}

\begin{center}
\begin{tabular}{@{}ccc@{}}

\includegraphics[width=22.5mm]{pictures/balanced/balancedunitleft1.pdf} & \raisebox{0.45\prelim}{$=$} &

\includegraphics[width=22.5mm]{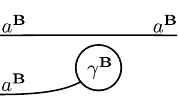}

\end{tabular}
\end{center}

holds in $Hom_{\mathfrak{A}}(\mathbf{B}(M\Box I\Box C, N), \mathbf{B}(M, C\Box N))$.

\end{enumerate}
\end{Definition}

\begin{Definition}
Let $\mathfrak{A}$ be a 2-category, and let $\mathbf{A}$, $\mathbf{B}$ be two $\mathfrak{C}$-balanced 2-functors $(\mathfrak{M},\mathfrak{N})\rightarrow \mathfrak{A}$. A $\mathfrak{C}$-balanced 2-natural transformation $\mathbf{A}\Rightarrow \mathbf{B}$ consists of:
\begin{enumerate}
    \item A 2-natural transformation $t:\mathbf{A}\Rightarrow \mathbf{B}$;
    \item An invertible modifications $\Pi^t$ given on $M$ in $\mathfrak{M}$, $C$ in $\mathfrak{C}$, and $N$ in $\mathfrak{N}$ by
\end{enumerate}

\begin{equation*}\begin{tikzcd}
{\mathbf{A}(MC,N)} \arrow[r, "a^{\mathbf{A}}"] \arrow[d, "t"']      & {\mathbf{A}(M,CN)} \arrow[d, "t"] \\
{\mathbf{B}(MC,N)} \arrow[r, "a^{\mathbf{B}}"'] \arrow[ru, "\Pi^t", shorten <= 2ex,shorten >= 2ex, Rightarrow] & {\mathbf{B}(M,CN);}               
\end{tikzcd}\end{equation*}

satisfying:

\begin{enumerate}
    \item[a.] For every $M$ in $\mathfrak{M}$, $C,D$ in $\mathfrak{C}$, and $N$ in $\mathfrak{N}$, the equality
\end{enumerate}

\settoheight{\prelim}{\includegraphics[width=52.5mm]{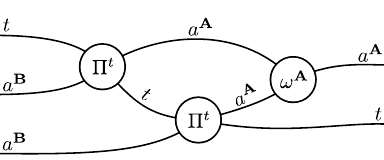}}

\begin{center}
\begin{tabular}{@{}ccc@{}}

\includegraphics[width=52.5mm]{pictures/balanced/balanced2nat1.pdf} & \raisebox{0.45\prelim}{$=$} &

\includegraphics[width=45mm]{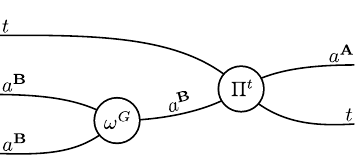}

\end{tabular}
\end{center}

\begin{enumerate}
    
    \item[] holds in $Hom_{\mathfrak{A}}(\mathbf{A}(M\Box C\Box D, N), \mathbf{B}(M, C\Box D\Box N))$,
    
    \item[b.] For every $M$ in $\mathfrak{M}$, and $N$ in $\mathfrak{N}$, the equality

\settoheight{\prelim}{\includegraphics[width=37.5mm]{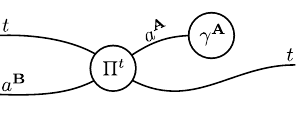}}

\begin{center}
\begin{tabular}{@{}ccc@{}}

\includegraphics[width=37.5mm]{pictures/balanced/balanced2nat3.pdf} & \raisebox{0.45\prelim}{$=$} &

\includegraphics[width=30mm]{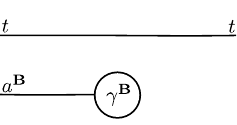}

\end{tabular}
\end{center}

holds in $Hom_{\mathfrak{A}}(\mathbf{A}(M\Box I, N), \mathbf{B}(M, N))$.

\end{enumerate}
\end{Definition}

\begin{Definition}
Let $\mathfrak{A}$ be a 2-category, $\mathbf{A}$, $\mathbf{B}$ be two $\mathfrak{C}$-balanced 2-functors $(\mathfrak{M},\mathfrak{N})\rightarrow \mathfrak{A}$, and $s$, $t$ be two $\mathfrak{C}$-balanced 2-natural transformations $\mathbf{A}\Rightarrow \mathbf{B}$. A $\mathfrak{C}$-balanced modification $s\Rrightarrow t$ is a modification $\delta:s\Rrightarrow t$ such that for every $M$ in $\mathfrak{M}$, $C$ in $\mathfrak{C}$, and $N$ in $\mathfrak{N}$ the equality

\settoheight{\prelim}{\includegraphics[width=30mm]{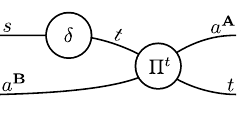}}

\begin{center}
\begin{tabular}{@{}ccc@{}}

\includegraphics[width=30mm]{pictures/balanced/balancedmodif1.pdf} & \raisebox{0.4\prelim}{$=$} & \includegraphics[width=30mm]{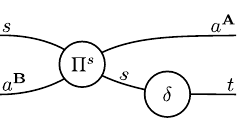}

\end{tabular}
\end{center}

\noindent holds in $Hom_{\mathfrak{A}}(\mathbf{A}(M\Box C,N), \mathbf{B}(M, C\Box N)).$
\end{Definition}

\subsection{Definition and Existence}

Let us now assume that $\mathfrak{C}$ is a compact semisimple tensor 2-category over a field $\mathds{k}$. Further, we will also fix $\mathfrak{M}$ a compact semisimple right $\mathfrak{C}$-module 2-category, and $\mathfrak{N}$ a compact semisimple left $\mathfrak{C}$-module 2-category. 

\begin{Definition}\label{def:relative2tensor}
If it exists, we write $\boxtimes_{\mathfrak{C}}: \mathfrak{M}\times\mathfrak{N}\rightarrow \mathfrak{M}\boxtimes_{\mathfrak{C}}\mathfrak{N}$ for the $\mathfrak{C}$-balanced 2-functor to a compact semisimple 2-category satisfying the following 3-universal property:
\begin{enumerate}
\item For any compact semisimple 2-category $\mathfrak{A}$, and every $\mathfrak{C}$-balanced bilinear 2-functor $\mathbf{B}:\mathfrak{M}\times\mathfrak{N}\rightarrow \mathfrak{A}$, there exists a linear 2-functor $\mathbf{B}':\mathfrak{M}\boxtimes_{\mathfrak{C}}\mathfrak{N}\rightarrow\mathfrak{A}$ and a $\mathfrak{C}$-balanced 2-natural equivalence $u:\mathbf{B}'\circ\boxtimes_{\mathfrak{C}}\Rightarrow \mathbf{B}$.

\item For every 2-functors $\mathbf{F},\mathbf{G}:\mathfrak{M}\boxtimes_{\mathfrak{C}}\mathfrak{N}\rightarrow \mathfrak{A}$, and $\mathfrak{C}$-balanced 2-natural transformation $t:\mathbf{F}\circ\boxtimes_{\mathfrak{C}}\Rightarrow \mathbf{G}\circ \boxtimes_{\mathfrak{C}}$, there exists a 2-natural transformation $t':\mathbf{F}\Rightarrow \mathbf{G}$ together with a $\mathfrak{C}$-balanced invertible modification $\eta:t' \circ \boxtimes_{\mathfrak{C}} \Rrightarrow t$.

\item Furthermore, for any 2-natural transformations $r,s:\mathbf{F}\Rightarrow \mathbf{G}$ and $\mathfrak{C}$-balanced modification $\delta:r\circ \boxtimes_{\mathfrak{C}}\Rrightarrow s\circ\boxtimes_{\mathfrak{C}}$, there exists a unique modification $\delta:r\Rrightarrow s$ such that $\delta'\circ\boxtimes_{\mathfrak{C}}=\delta$.
\end{enumerate}
\end{Definition}

\begin{Remark}
If we take $\mathfrak{C} = \mathbf{2Vect}_{\mathds{k}}$, in the above definition, we recover the definition of the 2-Deligne tensor product $\boxtimes$ from \cite{D3}, where it was denoted by $\boxdot$. Namely, there is an essentially unique $\mathbf{2Vect}_{\mathds{k}}$-module structure on any compact semisimple 2-category. The existence of the 2-Deligne tensor product was established in theorem 3.7 of the aforementioned reference.
\end{Remark}

\begin{Remark}\label{rem:2Delignereformulation}
Observe that $\mathfrak{C}$-balanced 2-functors, $\mathfrak{C}$-balanced 2-natural transformations, and $\mathfrak{C}$-balanced modifications can be assembled into a 2-category, which we denote by $\mathbf{Bal}_{\mathfrak{C}}(\mathfrak{M}\times\mathfrak{N},\mathfrak{A})$. Definition \ref{def:relative2tensor} may then be reformulated more concisely as saying that precomposition with $\boxtimes_{\mathfrak{C}}$ induces an equivalence of 2-categories \begin{equation*}\mathbf{Fun}(\mathfrak{M}\boxtimes_{\mathfrak{C}}\mathfrak{N},\mathfrak{A})\simeq \mathbf{Bal}_{\mathfrak{C}}(\mathfrak{M}\times\mathfrak{N},\mathfrak{A}),\end{equation*} from the 2-category of linear 2-functors $\mathfrak{M}\boxtimes_{\mathfrak{C}}\mathfrak{N}\rightarrow\mathfrak{A}$ to the 2-category of $\mathfrak{C}$-balanced bilinear 2-functors $\mathfrak{M}\times\mathfrak{N}\rightarrow\mathfrak{A}$. In fact, all the bilinear 2-functors that do occur may be factored through a 2-Deligne tensor product $\boxtimes$. As a consequence, the relative 2-Deligne tensor product may also be interpreted as a 3-colimit in the 3-category $\mathbf{CSS2C}$ of compact semisimple 2-categories \cite{D5}. More precisely, it is the geometric realization of the following diagram: \begin{equation*}\begin{tikzcd}[sep=small]
\mathfrak{MCCCN} \arrow[rr, shift left=6] \arrow[rr, shift left=2] \arrow[rr, shift right=2] \arrow[rr, shift right=6] &  & \mathfrak{MCCN} \arrow[rr, shift left=4] \arrow[rr, shift right=4] \arrow[rr] \arrow[ll, shift left=4] \arrow[ll, shift right=4] \arrow[ll] &  & \mathfrak{M\boxtimes C\boxtimes N} \arrow[rr, shift left=2] \arrow[rr, shift right=2] \arrow[ll, shift right=2] \arrow[ll, shift left=2] &  & \mathfrak{M\boxtimes N} \arrow[ll] \arrow[rr, dotted] &  & \mathfrak{M}\boxtimes_{\mathfrak{C}}\mathfrak{N}.
\end{tikzcd}\end{equation*}
\end{Remark}

We will now state our first main theorem, which is a relative version of the construction given in \cite{D3}. The decategorified version was studied in \cite{ENO2, DSPS14}.

\begin{Theorem}\label{thm:relative2Deligne}
For any separable compact semisimple right $\mathfrak{C}$-module 2-category $\mathfrak{M}$ and separable compact semisimple left $\mathfrak{C}$-module 2-category $\mathfrak{N}$, the relative 2-Deligne tensor product $\boxtimes_{\mathfrak{C}}: \mathfrak{M}\times\mathfrak{N}\rightarrow \mathfrak{M}\boxtimes_{\mathfrak{C}}\mathfrak{N}$ exists.
\end{Theorem}
\begin{proof}
By assumption, there exists separable algebras $A$ and $B$ in $\mathfrak{C}$ such that $\mathfrak{M}\simeq\mathbf{LMod}_{\mathfrak{C}}(A)$, the 2-category of left $A$-modules in $\mathfrak{C}$, and $\mathfrak{N}\simeq \mathbf{Mod}_{\mathfrak{C}}(B)$ as $\mathfrak{C}$-module 2-categories. Then, the bilinear 2-functor \begin{equation*}\begin{tabular}{r c c c}
$\mathbf{C}:$& $\mathfrak{M}\times\mathfrak{N}$ & $\rightarrow$ & $\mathbf{Bimod}_{\mathfrak{C}}(A,B)$\\
& $(M,N)$&$\mapsto$&$M\Box N$
\end{tabular}\end{equation*} is canonically $\mathfrak{C}$-balanced, and the target is a compact semisimple 2-category thanks to proposition 3.1.3 of \cite{D7}.

Now, let $\mathfrak{A}$ be any compact semisimple 2-category, and let $\mathbf{B}:(\mathfrak{M},\mathfrak{N})\rightarrow \mathfrak{A}$ be a $\mathfrak{C}$-balanced 2-functor. Because $\mathbf{Bimod}_{\mathfrak{C}}(A,B)$ is a compact semisimple 2-category, it is equivalent to the Cauchy completion of its full sub-2-category on the object $A\Box C\Box C\Box B$ for a well-chosen $C$ in $\mathfrak{C}$. For instance, one can take $C$ to be given by the direct sum of a simple object in each connected component of $\mathfrak{C}$. In particular, it follows from the 3-universal property of the Cauchy completion \cite{D1} that any linear 2-functor $\mathbf{B}':\mathbf{Bimod}_{\mathfrak{C}}(A,B)\rightarrow \mathfrak{A}$ is completely determined by its value on the full sub-2-category $\mathrm{B} End_{A\mathrm{-}B}(A\Box C\Box C\Box B)$ on the object $A\Box C\Box C\Box B$.

So let us turn our attention to defining the desired linear 2-functor \begin{equation*}\mathbf{B}':\mathrm{B} End_{A\mathrm{-}B}(A\Box C\Box C\Box B)\rightarrow \mathfrak{A}.\end{equation*} At the level of object, it is evident that we want $\mathbf{B}'(A\Box C\Box C\Box B)= \mathbf{B}(A\Box C, C\Box B)\simeq \mathbf{B}(A, C\Box C\Box B)$. However this assignment can not be extended directly to 1-morphisms:\ Not every $A$-$B$-bimodule 1-endomorphism of $A\Box C\Box C\Box B$ can be written as $f\Box g$ with $f:A\Box C\rightarrow A\Box C$ and $g:C\Box B\rightarrow C\Box B$. Instead, we will need to consider certain idempotents in the 1-category $End_{A\mathrm{-}B}(A\Box C\Box C\Box B)$. In order to define them, we will use the fact that $A$ and $B$ are separable algebras, that is, the counits $\epsilon_A:m_A\circ m^*_A\Rightarrow Id_A$ and $\epsilon_B:m_B\circ m^*_B\Rightarrow Id_B$, which we depict using cups, have sections $\gamma_A$ and $\gamma_B$ as bimodule 2-morphisms. We use the notations of \cite[section 2.1]{D7} for the associated coherence 2-isomorphisms.

Given any $A$-$B$-bimodule 1-morphism $h:A\Box C\Box C\Box B\rightarrow A\Box C\Box C\Box B$, we write $(e_h,\varepsilon_h)$ for the idempotent whose underlying $A$-$B$-bimodule 1-morphism is \begin{align*}e_h:ACCB\xrightarrow{ACCm^*_B}ACCBB&\xrightarrow{m^*_ACCm^*_BB}AACCBBB\xrightarrow{AhBB}AACCBBB\\ & \xrightarrow{m_ACCm_BB}ACCBB\xrightarrow{ACCm_B}ACCB,\end{align*} and with $A$-$B$-bimodule 2-morphism $\varepsilon_h$ depicted in figure \ref{fig:varepsilonh}. It is easy to check that $\varepsilon_h$ is an idempotent $A$-$B$-bimodule 2-morphism using the equations $\epsilon_A\circ \gamma_A = Id_{Id_A}$, $\epsilon_B\circ \gamma_B = Id_{Id_B}$, the facts that $m^*_A$ is a left $A$-module 1-morphism, $m^*_B$ is a right $B$-module 1-morphism, and $h$ is an $A$-$B$-bimodule 1-morphism (see definitions 2.1.1 and 1.3.5 of \cite{D7} for the equations expressed in our graphical calculus). One also checks that the idempotent $(e_h,\varepsilon_h)$ splits to $h$ in $End_{A\mathrm{-}B}(A\Box C\Box C\Box B)$ via the obvious $A$-$B$-bimodule 2-morphisms. It follows from its construction that we may evaluate the functor $\mathbf{B}$ on the idempotent $(e_h,\varepsilon_h)$, that is $\mathbf{B}(e_h,\varepsilon_h)$ is the idempotent in $\mathfrak{A}$ whose underlying 1-morphism $\mathbf{B}(AC,CB)\rightarrow \mathbf{B}(AC,CB)$ is \begin{align*}\mathbf{B}(AC,Cm_B) \circ a^{\mathbf{B}}\circ \mathbf{B}(m_ACCm_B,&B)\circ \mathbf{B}(AhB,B)\circ\\ &\circ \mathbf{B}(m^*_ACCm^*_B,B)\circ (a^{\mathbf{B}})^{\bullet}\circ\mathbf{B}(AC,Cm^*_B),\end{align*} and with 2-morphism depicted in figure \ref{fig:varepsilonhimage}. We note that it is not automatic that this 2-morphism is an idempotent. This can be checked using the same relations used to check that $\varepsilon_h$ is an idempotent together with the fact that $\mathbf{B}$ is $\mathfrak{C}$-balanced. The key observation being that all the relations that are used are either of the form $\alpha1=\beta1$ or $1\alpha=1\beta$ for some appropriate 2-morphisms $\alpha$ and $\beta$, and such relations can be evaluated under $\mathbf{B}$ (up to potentially inserting instances of $a^{\mathbf{B}}$ and $\omega^{\mathbf{B}}$).

Now, $\mathfrak{A}$ is locally idempotent complete, so that the idempotent $\mathbf{B}(e_h,\varepsilon_h)$ may be split to a 1-morphism \begin{equation*}\mathbf{B}'(h):\mathbf{B}(AC,CB)\rightarrow \mathbf{B}(AC,CB).\end{equation*} We will now argue that this defines a 2-functor $\mathbf{B}':\mathrm{B} End_{A\mathrm{-}B}(A\Box C\Box C\Box B)\rightarrow \mathfrak{A}$. In order to do so, it will be enough to prove that for any two $A$-$B$-bimodule 1-morphisms $h,k:A\Box C\Box C\Box B\rightarrow A\Box C\Box C\Box B$ the idempotent $\mathbf{B}(e_{h\circ k},\varepsilon_{h\circ k})$ refines the idempotent $\mathbf{B}(e_{h},\varepsilon_{h})\circ \mathbf{B}(e_{k},\varepsilon_{k})$. Namely, all the desired properties then follow readily from the uniqueness of splittings of idempotents. But, to prove our claim, it will suffice to show that the idempotent $(e_{h\circ k},\varepsilon_{h\circ k})$ refines $(e_{h},\varepsilon_{h})\circ (e_{k},\varepsilon_{k})$ in $\mathrm{B} End_{A\mathrm{-}B}(A\Box C\Box C\Box B)$ by $A$-$B$-bimodule 2-morphisms that can be evaluated under $\mathbf{B}$. More precisely, we want to find $A$-$B$-bimodule 2-morphisms $\varpi:e_h\circ e_k\Rightarrow e_{h\circ k}$ and $\varsigma:e_{h\circ k}\Rightarrow e_h\circ e_k$ such that \begin{equation*}\varpi\cdot\varsigma = Id_{e_{h\circ k}}\ \textrm{and}\ \varepsilon_{h}\circ \varepsilon_{k} = \varsigma\cdot \varepsilon_{h\circ k}\cdot \varpi.\end{equation*} These conditions are satisfied, for instance, by the $A$-$B$-bimodule 2-morphisms given in figures \ref{fig:varepsilonhsplitting1} and \ref{fig:varepsilonhsplitting2}.

\begin{landscape}
    \vspace*{\fill}
    \begin{figure}[!hbt]
    \centering
    \includegraphics[width=150mm]{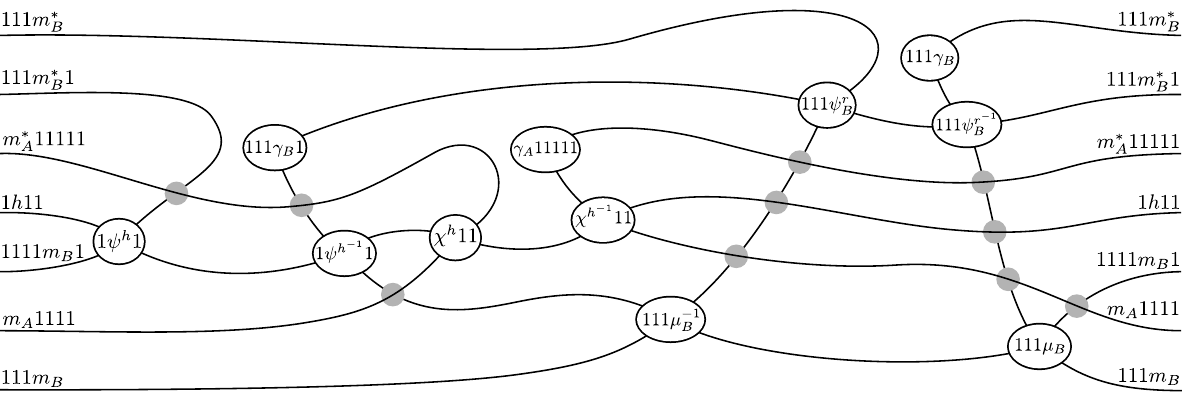}
    \caption{Definition of the idempotent $A$-$B$-bimodule 2-morphism $\varepsilon_h$}
    \label{fig:varepsilonh}
    \end{figure}
    \vfill
\end{landscape}

\begin{landscape}
    \vspace*{\fill}
    \begin{figure}[!hbt]
    \centering
    \includegraphics[width=165mm]{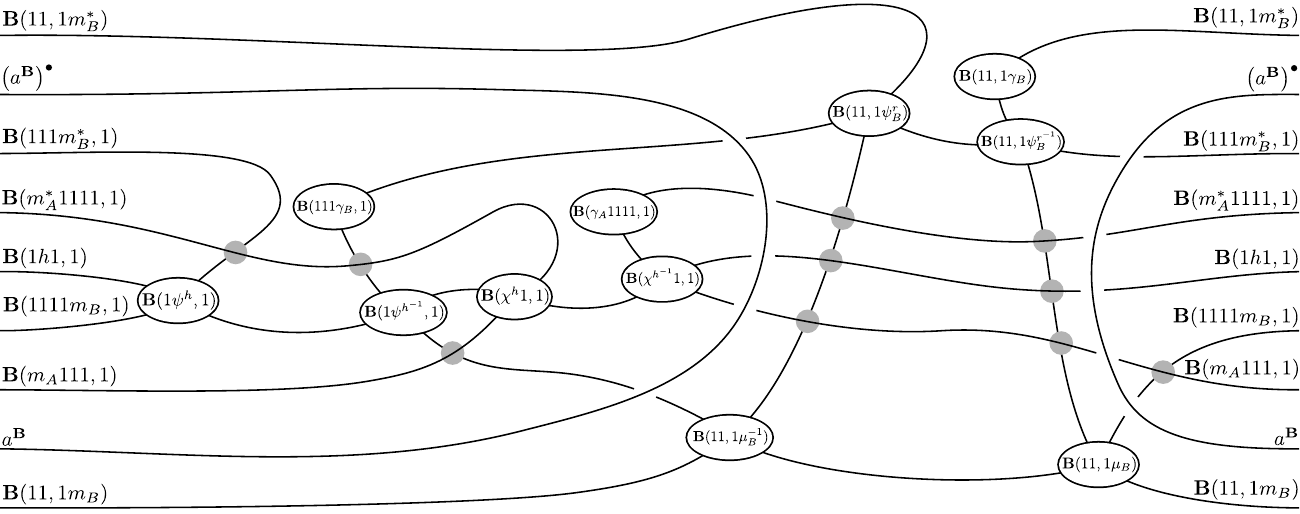}
    \caption{Image of the idempotent 2-morphism $\varepsilon_h$ under $\mathbf{B}'$}
    \label{fig:varepsilonhimage}
    \end{figure}
    \vfill
\end{landscape}

\begin{figure}[htb]
    \centering
    \includegraphics[width=120mm]{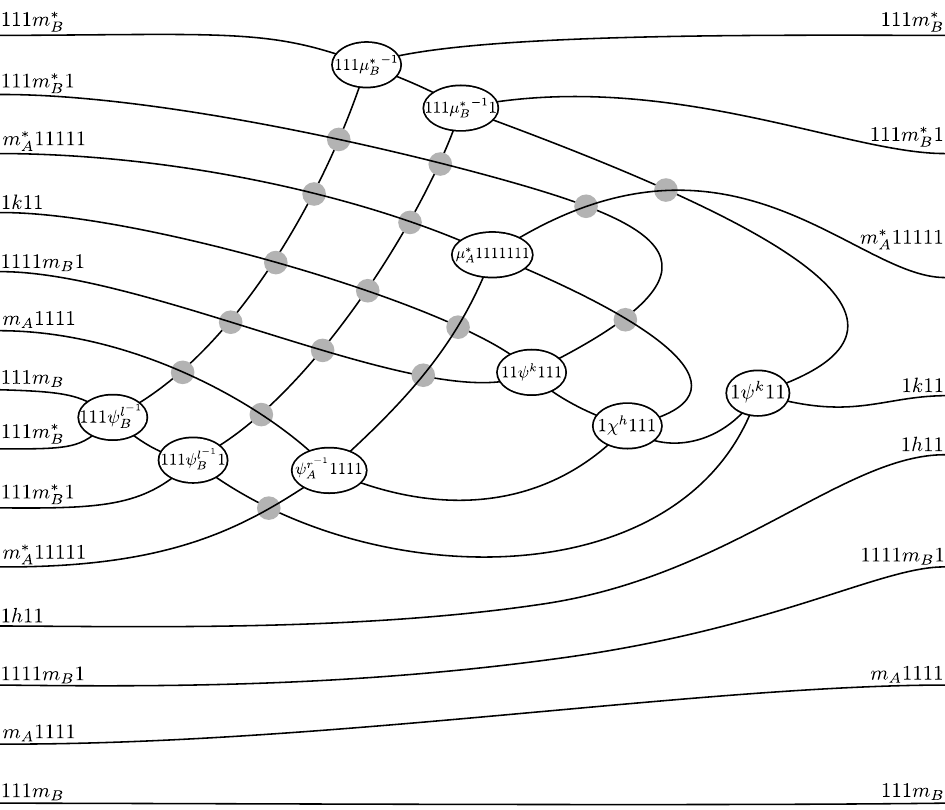}
    \caption{Refinement $A$-$B$-bimodule 2-morphism $\varpi:e_h\circ e_k\Rightarrow e_{h\circ k}$}
    \label{fig:varepsilonhsplitting1}
\end{figure}

\begin{figure}[hbt]
    \centering
    \includegraphics[width=120mm]{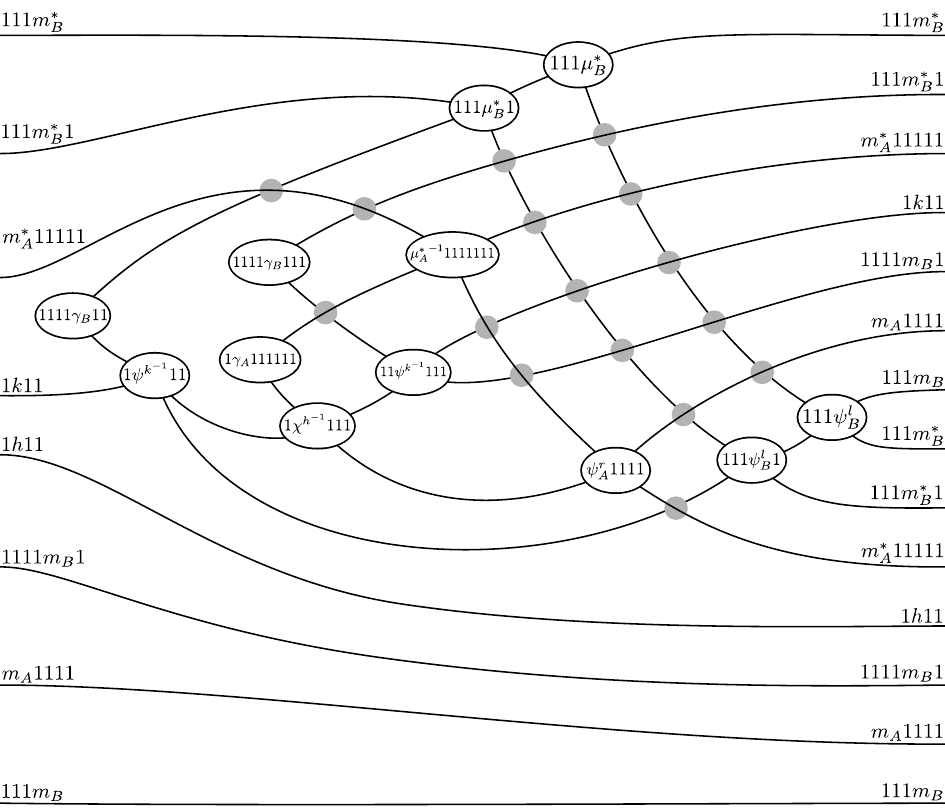}
    \caption{Refinement $A$-$B$-bimodule 2-morphism $\varsigma:e_{h\circ k}\Rightarrow e_h\circ e_k$}
    \label{fig:varepsilonhsplitting2}
\end{figure}

As $\mathbf{Bimod}_{\mathfrak{C}}(A,B)$ is the Cauchy completion of $\mathrm{B} End_{A\mathrm{-}B}(A\Box C\Box C\Box B)$, there is an essentially unique extension of $\mathbf{B}':\mathrm{B} End_{A\mathrm{-}B}(A\Box C\Box C\Box B)\rightarrow \mathfrak{A}$ to a linear 2-functor $\mathbf{B}':\mathbf{Bimod}_{\mathfrak{C}}(A,B)\rightarrow \mathfrak{A}$. Furthermore, it follows from the construction of $\mathbf{B}'$ that there is a $\mathfrak{C}$-balanced  2-natural equivalence $\mathbf{B}'\circ \mathbf{C}\simeq \mathbf{B}$. More precisely, let $M$ be any left $A$-module, and $N$ be any right $B$-module. By construction, $\mathbf{B}'\circ \mathbf{C}(M, N)$ is obtained by choosing a 2-condensation monad \begin{equation*}A\Box C\Box C\Box B\mathrel{\,\hspace{.75ex}\joinrel\rhook\joinrel\hspace{-.75ex}\joinrel\rightarrow} M\Box N\end{equation*} and then splitting in $\mathfrak{A}$ the image under $\mathbf{B}':\mathrm{B} End_{A\mathrm{-}B}(A\Box C\Box C\Box B)\rightarrow \mathfrak{A}$ of the underlying 2-condensation monad on $A\Box C\Box C\Box B$. Now, any two 2-condensation monads on $A\Box C\Box C\Box B$, which have the same splitting in $\mathbf{Bimod}_{\mathfrak{C}}(A,B)$ are necessarily equivalent by lemma \ref{lem:splittings2condensationmonads} above. But, thanks to our choice of $C$, we can always find 2-condensation monads $A\Box C\mathrel{\,\hspace{.75ex}\joinrel\rhook\joinrel\hspace{-.75ex}\joinrel\rightarrow} M$ and $C\Box B\mathrel{\,\hspace{.75ex}\joinrel\rhook\joinrel\hspace{-.75ex}\joinrel\rightarrow} N$. Taking the image under $\mathbf{C}$ of the product of these two 2-condensation monads, we obtain a 2-condensation monad \begin{equation*}A\Box C\Box C\Box B\mathrel{\,\hspace{.75ex}\joinrel\rhook\joinrel\hspace{-.75ex}\joinrel\rightarrow} M\Box N\end{equation*} of $A$-$B$-bimodules. But, taking the image under $\mathbf{B}$ of the product of these two 2-condensation monads, we obtain a 2-condensation monad \begin{equation*}\mathbf{B}(A\Box C,C\Box B)\mathrel{\,\hspace{.75ex}\joinrel\rhook\joinrel\hspace{-.75ex}\joinrel\rightarrow} \mathbf{B}(M,N)\end{equation*} in $\mathfrak{A}$. Putting the above discussion together, we find that there is canonical equivalence $\mathbf{B}'(M\Box N)\simeq \mathbf{B}(M,N)$ in $\mathfrak{A}$ as desired. A similar argument shows that these equivalences assemble to give a 2-natural equivalence and are compatible with the $\mathfrak{C}$-balanced structures.

Finally, it follows from variants of the above argument together with the 3-universal property of the Cauchy completion \cite{D1} that $\mathfrak{C}$-balanced 2-natural transformations and $\mathfrak{C}$-balanced modifications factor appropriately, which concludes the proof.
\end{proof}

\begin{Remark}
The proof of theorem \ref{thm:relative2Deligne} shows more generally that the $\mathfrak{C}$-balanced 2-functor $\boxtimes_{\mathfrak{C}}: \mathfrak{M}\times\mathfrak{N}\rightarrow \mathfrak{M}\boxtimes_{\mathfrak{C}}\mathfrak{N}$ is 3-universal with respect to any Cauchy complete 2-category $\mathfrak{A}$. Moreover, by inspecting the above proof very carefully, one finds that it is enough to assume that $\mathfrak{N}$ is a separable left $\mathfrak{C}$-module 2-category and that $\mathfrak{M}$ is a compact semisimple right $\mathfrak{C}$-module 2-category such that there exists a rigid algebra $A$ in $\mathfrak{C}$ such that $\mathfrak{M}\simeq\mathbf{LMod}_{\mathfrak{C}}(A)$ as right $\mathfrak{C}$-module 2-categories. (By theorem 5.3.4 and remark 5.3.10 of \cite{D4}, this last condition holds if $\mathds{k}$ is perfect and both $\mathfrak{C}$ and $\mathfrak{M}$ are locally separable.) Namely, using that $\mathfrak{M}\simeq\mathbf{LMod}_{\mathfrak{C}}(A)$ is compact semisimple by hypothesis, a slight generalization of proposition 3.1.2 of \cite{D7} establishes that $\mathbf{Bimod}_{\mathfrak{C}}(A,B)$ is compact semisimple. In addition, in the proof, one only ever uses that the counit $\epsilon_A:m_A\circ m_A^*\Rightarrow Id_A$ has a right adjoint as a left $A$-module, which also follows from the fact that $\mathbf{LMod}_{\mathfrak{C}}(A)$ is compact semisimple.
\end{Remark}

\subsection{Properties and Examples}

Let us fix $\mathfrak{C}$, $\mathfrak{D}$, and $\mathfrak{E}$ compact semisimple tensor 2-categories over the field $\mathds{k}$. Recall that a bimodule 2-category is \textit{separable} if it is separable independently as a left module 2-category and as a right module 2-category. We now derive some properties of separable compact semisimple bimodule 2-categories that are relevant to the proof of theorem \ref{thm:fulldualizability} below. Our discussion is inspired by \cite[section 2.4]{DSPS13}.

\begin{Proposition}\label{prop:separablebimodulecompose}
Let $\mathfrak{M}$ be a separable $\mathfrak{C}$-$\mathfrak{D}$-bimodule 2-category, and $\mathfrak{N}$ a separable $\mathfrak{D}$-$\mathfrak{E}$-bimodule 2-category. Then, the relative 2-Deligne tensor product $\mathfrak{M}\boxtimes_{\mathfrak{D}}\mathfrak{N}$ is a separable $\mathfrak{C}$-$\mathfrak{E}$-bimodule 2-category.
\end{Proposition}
\begin{proof}
We have seen in theorem \ref{thm:relative2Deligne} that $\mathfrak{M}\boxtimes_{\mathfrak{D}}\mathfrak{N}$ is a compact semisimple 2-category. Further, it follows readily from the 3-universal property that $\mathfrak{M}\boxtimes_{\mathfrak{D}}\mathfrak{N}$ inherits a $\mathfrak{C}$-$\mathfrak{E}$-bimodule structure, which is not a priori strict cubical. In any case, the $\mathfrak{D}$-balanced 2-functor $\boxtimes_{\mathfrak{D}}:\mathfrak{M}\times \mathfrak{N}\rightarrow \mathfrak{M}\boxtimes_{\mathfrak{D}}\mathfrak{N}$ is compatible with the $\mathfrak{C}$-$\mathfrak{E}$-bimodule structures.

Let $A$ be a separable algebra in $\mathfrak{C}$, and $B$ be a separable algebra in $\mathfrak{D}$ such that $\mathfrak{M}\simeq \mathbf{Mod}_{\mathfrak{C}}(A)$, and $\mathfrak{N}\simeq \mathbf{Mod}_{\mathfrak{D}}(B)$. It follows from the proof of theorem \ref{thm:relative2Deligne} that $\mathfrak{M}\boxtimes_{\mathfrak{D}}\mathfrak{N}\simeq \mathbf{Mod}_{\mathfrak{M}}(B)$. Moreover, the right action of $\mathfrak{D}$ on $\mathfrak{M}$ induces a monoidal 2-functor $F:\mathfrak{D}^{mop}\rightarrow \mathbf{End}_{\mathfrak{C}}(\mathfrak{M})\simeq \mathbf{Bimod}_{\mathfrak{C}}(A)$ by theorem 5.1.2 of \cite{D8}. In particular, we may think of the separable algebra $F(B)$ in $\mathbf{Bimod}_{\mathfrak{C}}(A)$ as a separable algebra in $\mathfrak{C}$ equipped with an algebra 1-homomorphism from $A$. Doing so, we find that $\mathbf{Mod}_{\mathfrak{M}}(B)\simeq \mathbf{Mod}_{\mathfrak{C}}(F(B))$ as left $\mathfrak{C}$-module 2-categories, so that $\mathfrak{M}\boxtimes_{\mathfrak{D}}\mathfrak{N}$ is separable as a left $\mathfrak{C}$-module 2-category. One checks similarly that $\mathfrak{M}\boxtimes_{\mathfrak{D}}\mathfrak{N}$ is separable as a right $\mathfrak{E}$-module 2-category.
\end{proof}

\begin{Example}
For simplicity, let us assume that $\mathds{k}$ is an algebraically closed field of characteristic zero. Let $\mathcal{B}$ be a braided fusion 1-category. We also consider $\mathcal{C}$ a $\mathcal{B}^{rev}$-central fusion 1-category, that is a fusion 1-category $\mathcal{C}$ equipped with a braided functor $\mathcal{B}^{rev}\rightarrow \mathcal{Z}(\mathcal{C})$ to its Drinfeld center, as well as $\mathcal{D}$ a $\mathcal{B}$-central fusion 1-category. In particular, $\mathbf{Mod}(\mathcal{B})$ is a fusion 2-category, $\mathbf{Mod}(\mathcal{C})$ is a finite semisimple right $\mathbf{Mod}(\mathcal{B})$-module 2-category, and $\mathbf{Mod}(\mathcal{D})$ is a finite semisimple left $\mathbf{Mod}(\mathcal{B})$-module 2-category. Then, it follows from the proof of theorem \ref{thm:relative2Deligne} that \begin{equation*}\mathbf{Mod}(\mathcal{C})\boxtimes_{\mathbf{Mod}(\mathcal{B})}\mathbf{Mod}(\mathcal{D})\simeq \mathbf{Mod}(\mathcal{C}\boxtimes_{\mathcal{B}}\mathcal{D})\end{equation*} as finite semisimple 2-categories. Namely, $\mathcal{C}\boxtimes_{\mathcal{B}}\mathcal{D}$ inherits a canonical structure of multifusion 1-category from the central structures on $\mathcal{C}$ and $\mathcal{D}$ (see for instance \cite{BJS}). Further, the 2-category of $\mathcal{C}$-$\mathcal{D}$-bimodules in $\mathbf{Mod}(\mathcal{B})$ is equivalent to the 2-category of finite semisimple right $\mathcal{C}\boxtimes_{\mathcal{B}}\mathcal{D}$-module 1-categories (see also theorem 3.27 of \cite{Lau} for a closely related observation). We note that this finite semisimple 2-category is in general not connected, i.e.\ $\mathcal{C}\boxtimes_{\mathcal{B}}\mathcal{D}$ is in general not a fusion 1-category. For instance, if $\mathcal{B}=\mathbf{Rep}(G)$ for some finite group $G$, and $\mathcal{C} = \mathcal{D} = \mathbf{Vect}$, then $\mathbf{Vect}\boxtimes_{\mathbf{Rep}(G)}\mathbf{Vect}=\boxplus_{g\in G}\mathbf{Vect}$.
\end{Example}

Recall that $\mathfrak{C}$ is in particular a rigid monoidal 2-category. This means that every object $C$ has a coherent right dual $(C,C^{\sharp},i_C,e_C,E_C,F_C)$. In fact, such a right dual is essentially unique \cite{Pstr}. It was shown in \cite[Appendix A.2]{D2} and \cite[Appendix A.1]{DX} that the assignment sending an object $C$ in $\mathfrak{C}$ to its right dual $C^{\sharp}$ defines a monoidal 2-functor $(-)^{\sharp}:\mathfrak{C}\rightarrow \mathfrak{C}^{1op,mop}$ to the monoidal 2-category obtained from $\mathfrak{C}$ by reversing the order of the monoidal product and the direction of 1-morphisms. Likewise, there is also a monoidal 2-functor ${^{\sharp}(-)}$ sending an object to its left dual.

\begin{Construction}
Let $\mathfrak{M}$ be a separable $\mathfrak{C}$-$\mathfrak{D}$-bimodule 2-category. We define a $\mathfrak{D}$-$\mathfrak{C}$-bimodule 2-category $\mathfrak{M}^*$ whose underlying 2-category is $\mathfrak{M}^{1op}$, the 2-category obtained from $\mathfrak{M}$ by reversing the direction of 1-morphisms. Given an object $M$ in $\mathfrak{M}$, we use $M^*$ to denote the corresponding object of $\mathfrak{M}^*$. Then, the $\mathfrak{D}$-$\mathfrak{C}$-bimodule structure on $\mathfrak{M}^*$ is given by \begin{equation*}D\Box M^*\Box C := ({^{\sharp}C}\Box M\Box {^{\sharp}D})^*\end{equation*} for any objects $M^*$ in $\mathfrak{M}^*$, $C$ in $\mathfrak{C}$, and $D$ in $\mathfrak{D}$. Likewise, we define the $\mathfrak{D}$-$\mathfrak{C}$-bimodule 2-category ${^*\mathfrak{M}}$ as $\mathfrak{M}^{1op}$ with bimodule structure given by \begin{equation*}D\Box{^*M}\Box C := {^*({C^{\sharp}}\Box M\Box {D^{\sharp}})}\end{equation*} for any objects ${^*M}$ in $^*\mathfrak{M}$, $C$ in $\mathfrak{C}$, and $D$ in $\mathfrak{D}$.
\end{Construction}

\begin{Lemma}\label{lem:functorduals}
Assume that $\mathds{k}$ is perfect. Let $\mathfrak{C}$ and $\mathfrak{D}$ be a locally separable compact semisimple tensor 2-categories, and let $\mathfrak{M}$ be a separable $\mathfrak{C}$-$\mathfrak{D}$-bimodule 2-category. Then we have \begin{equation*}\mathbf{Fun}_{\mathfrak{C}}(\mathfrak{M},\mathfrak{C})\simeq {^*\mathfrak{M}}\ \textrm{and}\ \mathbf{Fun}_{\mathfrak{D}}(\mathfrak{M},\mathfrak{D})\simeq {\mathfrak{M}^*}\end{equation*} as $\mathfrak{D}$-$\mathfrak{C}$-bimodule 2-categories.
\end{Lemma}
\begin{proof}
Consider the 2-functor \begin{equation*}\begin{tabular}{r c c c}
$\mathbf{A}:$&$\mathbf{Fun}_{\mathfrak{C}}(\mathfrak{M},\mathfrak{C})$&$\rightarrow$&$\mathfrak{M}^{1op}$\\ & $F$&$\mapsto$&$F^*(I)$
\end{tabular}\end{equation*} given by sending a left $\mathfrak{C}$-module 2-functor $F:\mathfrak{M}\rightarrow \mathfrak{C}$ to its right adjoint $F^*$ evaluated at the identity object $I$ in $\mathfrak{C}$. This right adjoint exists thanks to corollary 4.2.3 of \cite{D5}, and it inherits a canonical left $\mathfrak{C}$-module structure by proposition 4.2.3 of \cite{D8}. This 2-functor is an equivalence. Namely, it is the composite of the 2-functor $(-)^*:\mathbf{Fun}_{\mathfrak{C}}(\mathfrak{M},\mathfrak{C})\rightarrow \mathbf{Fun}_{\mathfrak{C}}(\mathfrak{C},\mathfrak{M})^{1op}$ sending a functor to its right adjoint, and the evaluation 2-functor $\mathbf{E}:\mathbf{Fun}_{\mathfrak{C}}(\mathfrak{C},\mathfrak{M})^{1op}\rightarrow \mathfrak{M}^{1op}$;\ Both of which are manifestly equivalences.

It therefore only remains to check that $\mathbf{A}$ is compatible with the bimodule structures. The right $\mathfrak{C}$-action on $\mathbf{Fun}_{\mathfrak{C}}(\mathfrak{M},\mathfrak{C})$ is given by postcomposition with $(-)\Box C$. But the right adjoint of $(-)\Box C$ is $(-)\Box C^{\sharp}$. Likewise, the left $\mathfrak{D}$-action on $\mathbf{Fun}_{\mathfrak{C}}(\mathfrak{M},\mathfrak{C})$ is given by precomposition with $(-)\Box D$, and the right adjoint of $(-)\Box D$ is again $(-)\Box D^{\sharp}$. This shows that the 2-functor $\mathbf{A}$ admits bimodule structure as desired. The other equivalence is derived analogously.
\end{proof}

\begin{Lemma}\label{lem:moduleduals}
Let $A$ be a separable algebra in a compact semisimple tensor 2-category $\mathfrak{C}$. Then, we have equivalences of $\mathfrak{C}$-module 2-categories \begin{equation*}\mathbf{Mod}_{\mathfrak{C}}(A)^*\simeq \mathbf{LMod}_{\mathfrak{C}}(A^{\sharp\sharp}),\end{equation*}
\begin{equation*}\mathbf{LMod}_{\mathfrak{C}}(A)^*\simeq \mathbf{Mod}_{\mathfrak{C}}(A)\end{equation*} given by sending an $A$-module $M^*$ to $M^{\sharp}$, and 
\begin{equation*}{^*\mathbf{Mod}_{\mathfrak{C}}(A)}\simeq \mathbf{LMod}_{\mathfrak{C}}(A),\end{equation*} \begin{equation*}{^*\mathbf{LMod}_{\mathfrak{C}}(A)}\simeq \mathbf{Mod}_{\mathfrak{C}}({^{\sharp\sharp}A}),\end{equation*} given by sending an $A$-module $^*M$ to $^{\sharp}M$.
\end{Lemma}
\begin{proof}
We prove the first equivalence, the others follow similarly. Let $M$ be a right $A$-module. We write $M^*$ for the corresponding object of $\mathbf{Mod}_{\mathfrak{C}}(A)^*$. Then, applying $(-)^{\sharp}$ to the action 1-morphism $n^M:M\Box A\rightarrow M$ yields a 1-morphism $(n^M)^{\sharp}:M^{\sharp}\rightarrow A^{\sharp}\Box M^{\sharp}$. Upon taking its mate, we get a 1-morphism \begin{equation*}l^{M^{\sharp}}:A^{\sharp\sharp}\Box M^{\sharp}\xrightarrow{A^{\sharp\sharp}\Box(n^M)^{\sharp}}A^{\sharp\sharp}\Box A^{\sharp}\Box M^{\sharp}\xrightarrow{(i_{A^{\sharp}})^*\Box M^{\sharp}}M^{\sharp}.\end{equation*} Similarly, the coherence 2-isomorphisms for the right $A$-action on $M$ yield coherence 2-isomorphisms for a left $A^{\sharp\sharp}$-action on $M^{\sharp}$. Further, this construction can be extended in the obvious way to 1-morphisms (remembering that the underlying 2-category of $\mathbf{Mod}_{\mathfrak{C}}(A)^*$ is $\mathbf{Mod}_{\mathfrak{C}}(A)^{1op}$), and is functorial. As the monoidal 2-functor $(-)^{\sharp}:\mathfrak{C}\rightarrow \mathfrak{C}^{mop,1op}$ is an equivalence, so is the 2-functor $M^*\mapsto M^{\sharp}$. Further, it is clear that this is suitably compatible with the left $\mathfrak{C}$-module structures, which concludes the proof of the lemma.
\end{proof}

As already used above, the property of being separable for an algebra in a monoidal 2-category is preserved by all monoidal 2-functors. In particular, if $A$ is a separable algebra in a compact semisimple tensor 2-category $\mathfrak{C}$, then $A^{\sharp\sharp}$ is separable. Thus, we obtain the following corollary of the two lemmas above.

\begin{Corollary}\label{cor:adjointbimoduleseparable}
Let $\mathfrak{C}$ and $\mathfrak{D}$ be locally separable compact semisimple tensor 2-categories over the perfect field $\mathds{k}$. For any separable $\mathfrak{C}$-$\mathfrak{D}$-bimodule 2-category $\mathfrak{M}$, we have that $\mathbf{Fun}_{\mathfrak{C}}(\mathfrak{M},\mathfrak{C})$ is a separable $\mathfrak{D}$-$\mathfrak{C}$-bimodule 2-category.
\end{Corollary}

\begin{Proposition}\label{prop:modulefunctorequivalence}
Assume that $\mathds{k}$ is perfect. Let $\mathfrak{C}$, $\mathfrak{D}$, and $\mathfrak{E}$ be locally separable compact semisimple tensor 2-categories. Further, let $\mathfrak{M}$ be a separable $\mathfrak{C}$-$\mathfrak{D}$-bimodule 2-category, and $\mathfrak{N}$ a separable $\mathfrak{C}$-$\mathfrak{E}$-bimodule 2-category. Then, the canonical $\mathfrak{C}$-balanced 2-functor $\mathbf{Fun}_{\mathfrak{C}}(\mathfrak{M},\mathfrak{C})\times\mathfrak{N}\rightarrow \mathbf{Fun}_{\mathfrak{C}}(\mathfrak{M},\mathfrak{N})$ given by $(F,N)\mapsto F(-)\Box N$ induces an equivalence of $\mathfrak{D}$-$\mathfrak{E}$-bimodule 2-categories, \begin{equation*}\mathbf{Fun}_{\mathfrak{C}}(\mathfrak{M},\mathfrak{C})\boxtimes_{\mathfrak{C}}\mathfrak{N}\rightarrow \mathbf{Fun}_{\mathfrak{C}}(\mathfrak{M},\mathfrak{N}).\end{equation*}
\end{Proposition}
\begin{proof}
Let $A$ and $B$ be separable algebras in $\mathfrak{C}$ such that $\mathfrak{M}\simeq\mathbf{Mod}_{\mathfrak{C}}(A)$ and $\mathfrak{N}\simeq\mathbf{Mod}_{\mathfrak{C}}(B)$ as left $\mathfrak{C}$-module 2-categories. It follows from theorem 5.1.2 of \cite{D8} that $\mathbf{Fun}_{\mathfrak{C}}(\mathfrak{M},\mathfrak{C})\simeq \mathbf{LMod}_{\mathfrak{C}}(A)$ as right $\mathfrak{C}$-module 2-categories, and that there is an equivalence of plain 2-categories $\mathbf{Fun}_{\mathfrak{C}}(\mathfrak{M},\mathfrak{N})\simeq \mathbf{Bimod}_{\mathfrak{C}}(A,B)$. The equivalence of 2-categories in the statement of the proposition then follows from the construction given in theorem \ref{thm:relative2Deligne}, and the fact that it is compatible with the $\mathfrak{D}$-$\mathfrak{E}$-bimodule structures follows from the 3-universal property of the relative 2-Deligne tensor product.
\end{proof}

\section{Separability and Full Dualizability}

\subsection{Separability}

We begin by discussing the generalization to compact semisimple tensor 2-categories over a perfect field $\mathds{k}$ of the notion of separability introduced in \cite{D9} for fusion 2-category over an algebraically closed field of characteristic zero. This categorifies definition 2.5.8 of \cite{DSPS13}.

\begin{Definition}
A compact semisimple tensor 2-category $\mathfrak{C}$ is called separable if it is locally separable and it is separable as a left $\mathfrak{C}\boxtimes\mathfrak{C}^{mop}$-module 2-category.
\end{Definition}

Categorifying corollary 2.5.9 of \cite{DSPS13}, we have the following characterization of separability. The case of algebraically closed fields of characteristic zero was given in construction 5.2.2 of \cite{D9}.

\begin{Proposition}
A locally separable compact semisimple tensor 2-category is separable if and only if its Drinfeld center is compact semisimple.
\end{Proposition}
\begin{proof}
By lemma 2.2.1 of \cite{D9}, the result is a direct consequence of proposition 5.1.5 of \cite{D8}.
\end{proof}

\begin{Remark}
Over fields of positive characteristic, there exists compact semisimple tensor 2-categories that are not separable. For instance, if $\mathds{k}$ is an algebraically closed field of characteristic $p>0$, then the locally separable compact semisimple tensor 2-category $\mathbf{2Vect}(\mathbb{Z}/p)$ of $\mathbb{Z}/p$-graded 2-vector spaces is not separable. Namely, over algebraically closed fields of characteristic zero a notion of dimension for fusion 2-categories, whose non-vanishing is equivalent to separability, was introduced in \cite{D9} using the notion of dimension of a connected rigid algebra from \cite[section 3.2]{D7}. It is straightforward to generalize these considerations to any algebraically closed field. Further, we have $\mathrm{Dim}(\mathbf{2Vect}(\mathbb{Z}/p))=p=0$ in $\mathds{k}$, so that $\mathbf{2Vect}(\mathbb{Z}/p)$ is not separable. In fact, it can also be checked directly that there is an equivalence of 2-categories \begin{equation*}\mathscr{Z}(\mathbf{2Vect}(\mathbb{Z}/p))\simeq \boxplus_{i=1}^p\mathbf{2Rep}(\mathbb{Z}/p),\end{equation*} where $\mathbf{2Rep}(\mathbb{Z}/p)$ is the 2-category of 2-representations of $\mathbb{Z}/p$ in $\mathbf{2Vect}$. In particular, as $End_{\mathbf{2Rep}(\mathbb{Z}/p)}(\mathbf{Vect})\simeq \mathbf{Rep}(\mathbb{Z}/p)$ is not semisimple, $\mathscr{Z}(\mathbf{2Vect}(\mathbb{Z}/p))$ is not semisimple as expected.
\end{Remark}

The following result is an internalization of corollary 2.6.8 of \cite{DSPS13}, and generalizes theorem 5.1.1 of \cite{D9}.

\begin{Theorem}
Over a field of characteristic zero, every rigid algebra in a compact semisimple tensor 2-category is separable.
\end{Theorem}
\begin{proof}
Let $\mathfrak{C}$ be a compact semisimple tensor 2-category over a field $\mathds{k}$ of characteristic zero. Let $\overline{\mathds{k}}$ be the algebraic closure of $\mathds{k}$. We write $\mathfrak{C}\widehat{\otimes}\overline{\mathds{k}}$ for the local Cauchy completion of the linear tensor product $\mathfrak{C}\otimes\overline{\mathds{k}}$. Then, we set $\overline{\mathfrak{C}}:= Cau(\mathfrak{C}\widehat{\otimes}\overline{\mathds{k}})$, which is a compact semisimple tensor 2-category over $\overline{\mathds{k}}$ (see \cite{D3}). Let us consider a rigid algebra $A$ in $\mathfrak{C}$. We claim that there is a canonical equivalence \begin{equation*}\overline{\mathbf{Bimod}_{\mathfrak{C}}(A)}:=Cau(\mathbf{Bimod}_{\mathfrak{C}}(A)\widehat{\otimes}\overline{\mathds{k}})\simeq \mathbf{Bimod}_{\overline{\mathfrak{C}}}(A)\end{equation*} of 2-categories. Then, let $Z(A)$ denote the $\mathds{k}$-linear 1-category of $A$-$A$-bimodule endomorphisms of $A$ in $\mathfrak{C}$, which is Cauchy complete and locally finite. Assuming that the above equivalence holds, we have that $Cau(Z(A)\otimes\overline{\mathds{k}})$ is a finite semisimple 1-category by theorem 5.1.1 of \cite{D9}. As $\mathds{k}$ is perfect, this implies that $Z(A)$ is finite semisimple. It then follows from theorem 3.1.6 of \cite{D7} that $A$ is separable.

It remains to prove the claim. We begin by observing that the canonical $\overline{\mathds{k}}$-linear 2-functor $\mathbf{F}:\mathbf{Bimod}_{\mathfrak{C}}(A)\widehat{\otimes}\overline{\mathds{k}}\rightarrow \mathbf{Bimod}_{\overline{\mathfrak{C}}}(A)$ is fully faithful on 2-morphisms. In fact, as being fully faithful on 2-morphisms is preserved by local Cauchy completion, it is enough to prove that the canonical $\overline{\mathds{k}}$-linear 2-functor $\mathbf{Bimod}_{\mathfrak{C}}(A)\otimes\overline{\mathds{k}}\rightarrow \mathbf{Bimod}_{\overline{\mathfrak{C}}}(A)$ is fully faithful on 2-morphisms. To see this, let $f,g$ be two $A$-$A$-bimodule 1-morphisms in $\mathfrak{C}$. By definition, the $\mathds{k}$-vector space of $A$-$A$-bimodule 2-morphisms $f\Rightarrow g$ in $\mathfrak{C}$ is a linear subspace of $End_{\mathfrak{C}}(f,g)$ cut out by certain linear equations. But, the $\overline{\mathds{k}}$-vector space of $A$-$A$-bimodule 2-morphisms $f\Rightarrow g$ in $\overline{\mathfrak{C}}$ is the linear subspace of $End_{\overline{\mathfrak{C}}}(f,g)=End_{\mathfrak{C}}(f,g)\otimes \overline{\mathds{k}}$ cut out by precisely the same linear equations. This shows that the 2-functor under consideration is indeed fully faithful on 2-morphisms.

We now argue that $\mathbf{F}$ is a Cauchy completion. In order to see this, it is enough to show that $\mathbf{Bimod}_{\mathfrak{C}}(A)\widehat{\otimes}\overline{\mathds{k}}$ contains a full sub-2-category, whose Cauchy completion is $\mathbf{Bimod}_{\overline{\mathfrak{C}}}(A)$. To this end, we will consider $\mathfrak{F}$, the sub-2-category of $\mathbf{Bimod}_{\mathfrak{C}}(A)\widehat{\otimes}\overline{\mathds{k}}$ on the free objects, i.e.\ the bimodules of the form $A\Box C\Box A$ for some $C$ in $\mathfrak{C}$, which is full on both 1-morphisms and 2-morphisms. Given any objects $C$ of $\mathfrak{C}$ and $P$ of $\mathbf{Bimod}_{\mathfrak{C}}(A)\widehat{\otimes}\overline{\mathds{k}}$, the functor \begin{equation*}Cau(Hom_{A\textrm{-}A}(A\Box C\Box A, P)\otimes\overline{\mathds{k}})\rightarrow Hom_{A\textrm{-}A}(\mathbf{F}(A\Box C\Box A), \mathbf{F}(P))\end{equation*} induced by $\mathbf{F}$ is identified with the functor \begin{equation*}Cau(Hom_{\mathfrak{C}}(C, P)\otimes\overline{\mathds{k}})\rightarrow Hom_{\overline{\mathfrak{C}}}(C, P),\end{equation*} which is an equivalence of 1-categories by definition. Thus, the canonical inclusion $\mathfrak{F}\hookrightarrow \mathbf{Bimod}_{\overline{\mathfrak{C}}}(A)$ is full on 1-morphisms. Finally, the bimodule version of the proof of proposition 3.1.2 of \cite{D7} shows that the finite semisimple 2-category $\mathbf{Bimod}_{\overline{\mathfrak{C}}}(A)$ is generated under Cauchy completion by its full sub-2-category on the free objects. But, the canonical 2-functor $\mathfrak{C}\widehat{\otimes}\overline{\mathds{k}}\rightarrow\overline{\mathfrak{C}}$ is by definition a Cauchy completion, so that every free bimodule in $\mathbf{Bimod}_{\overline{\mathfrak{C}}}(A)$ is contained in the Cauchy completion of $\mathfrak{F}$. This shows that the Cauchy completion of $\mathfrak{F}$ is $\mathbf{Bimod}_{\overline{\mathfrak{C}}}(A)$ as desired, and therefore finishes the proof of the theorem.
\end{proof}

By theorem 5.3.2 of \cite{D8} and remark 5.3.10 of \cite{D4}, we then obtain:

\begin{Corollary}
Over a field of characteristic zero, for any compact semisimple tensor 2-category $\mathfrak{C}$, and compact semisimple left $\mathfrak{C}$-module 2-category $\mathfrak{M}$, we have that $\mathbf{End}_{\mathfrak{C}}(\mathfrak{M})$ is a compact semisimple tensor 2-category.
\end{Corollary}

The next result follows immediately.

\begin{Corollary}\label{cor:separabilitycharzero}
Over a field of characteristic zero, every compact semisimple tensor 2-category is separable.
\end{Corollary}

\subsection{Full Dualizability}

From now on, we will always assume that $\mathds{k}$ is a perfect field. We write $\mathbf{CSS2C}$ for the 3-category of compact semisimple 2-categories. The 2-Deligne tensor product from \cite{D4} endows this 3-category with a symmetric monoidal structure thanks to its universal property. Alternatively, it was shown in \cite{D3} that $\mathbf{CSS2C}$ is equivalent to the 3-category of finite semisimple tensor categories and \textit{right separable} bimodules 1-categories $\mathbf{TC}^{rsep}$, which is a close relative of the 3-categories introduced in \cite[section 2.5]{DSPS13}. Further, it follows from the arguments therein that the Deligne tensor product endows this 3-category with a symmetric monoidal structure, so that $\mathbf{CSS2C}$ is a symmetric monoidal 3-category. In fact, it follows from the 3-universal property of the 2-Deligne tensor product that these symmetric monoidal structures do agree (up to canonical equivalences).

Thanks to theorem \ref{thm:relative2Deligne} above, the relative 2-Deligne tensor product of any two separable bimodule 2-categories exists. Moreover, it follows from proposition \ref{prop:separablebimodulecompose} that the relative 2-Deligne tensor product of any two separable bimodule 2-categories is again separable. Let us also note that the 2-Deligne tensor product commutes with the relative 2-Deligne tensor product. This follows readily from the explicit construction of the relative 2-Deligne tensor product combined with a slight generalization of lemma 1.4.4 of \cite{D9}. This argument also shows that the 2-Deligne tensor product of any two separable bimodule 2-categories is a separable bimodule 2-category. Thus, we can apply the construction of \cite{JFS}, so as to obtain a symmetric monoidal 4-category $\mathbf{CST2C}$ of compact semisimple tensor 2-categories, and separable bimodule 2-categories. We are more specifically interested in the symmetric monoidal full sub-4-category $\mathbf{CST2C}^{sep}$ on the separable compact semisimple tensor 2-categories.

The next result is a categorification of theorem 3.4.3 of \cite{DSPS13}, by which our proof is inspired. In fact, the result below subsumes the content of their theorem as $\Omega \mathbf{CST2C} = \mathbf{CSS2C}^{sep}$, the symmetric monoidal 3-category of locally separable compact semisimple 2-categories, which can be identified with $\mathbf{TC}^{sep}$, the symmetric monoidal 3-category of separable tensor 1-categories. It also positively answers a question raised in \cite{DR}.

\begin{Theorem}\label{thm:fulldualizability}
The symmetric monoidal 4-category $\mathbf{CST2C}^{sep}$ is fully dualizable.
\end{Theorem}
\begin{proof}
Let $\mathfrak{C}$ be a separable compact semisimple tensor 2-category. Then, $\mathfrak{C}$ is a separable left $\mathfrak{C}\boxtimes\mathfrak{C}^{mop}$-module 2-category, and also a separable right $\mathfrak{C}^{mop}\boxtimes\mathfrak{C}$-module 2-category. Further, as $\mathfrak{C}$ is locally separable, $\mathfrak{C}$ is also separable as a module 2-category over $\mathbf{2Vect}_{\mathds{k}}$. These 1-morphisms provide evaluation and coevaluation 1-morphisms witnessing that $\mathfrak{C}$ and $\mathfrak{C}^{mop}$ are dual. (In fact, that every object is 1-dualizable is a general phenomenon in higher Morita categories \cite{GwS}. Let us however emphasize that some care has to be taken because not every compact semisimple bimodule 2-category yields a 1-morphism in $\mathbf{CST2C}$.)

We now argue that every separable $\mathfrak{C}$-$\mathfrak{D}$-bimodule 2-category $\mathfrak{M}$ between locally separable compact semisimple tensor 2-categories $\mathfrak{C}$ and $\mathfrak{D}$ viewed as a 1-morphism of $\mathbf{CST2C}$ has a left adjoint. We will prove that the $\mathfrak{D}$-$\mathfrak{C}$-bimodule 2-category $\mathbf{Fun}_{\mathfrak{C}}(\mathfrak{M},\mathfrak{C})$ is a left adjoint for $\mathfrak{M}$ in $\mathbf{CST2C}$. Firstly, that $\mathbf{Fun}_{\mathfrak{C}}(\mathfrak{M},\mathfrak{C})$ is a separable bimodule 2-category was recorded in corollary \ref{cor:adjointbimoduleseparable}. Secondly, we have seen in proposition \ref{prop:modulefunctorequivalence} that \begin{equation*}\mathbf{Fun}_{\mathfrak{C}}(\mathfrak{M},\mathfrak{C})\boxtimes_{\mathfrak{C}}\mathfrak{M}\simeq \mathbf{End}_{\mathfrak{C}}(\mathfrak{M})\end{equation*} as $\mathfrak{D}$-$\mathfrak{D}$-bimodule 2-categories, where the left $\mathfrak{D}$-action comes from precomposition. Thirdly, write \begin{equation*}\begin{tabular}{r c c c} $\eta:$& $\mathfrak{D}$ & $\rightarrow$ & $\mathbf{End}_{\mathfrak{C}}(\mathfrak{M})$,\\ & $D$&$\mapsto$&$(-)\Box D$\end{tabular}\end{equation*} and note that this is a $\mathfrak{D}$-$\mathfrak{D}$-bimodule 2-functor. Further, there is also a $\mathfrak{C}$-$\mathfrak{C}$-bimodule 2-functor \begin{equation*}\begin{tabular}{r c c c} $\epsilon:$& $\mathfrak{M}\boxtimes_{\mathfrak{D}}\mathbf{Fun}_{\mathfrak{C}}(\mathfrak{M},\mathfrak{C})$ & $\rightarrow$ & $\mathfrak{C}$\end{tabular}\end{equation*} induced by the canonical $\mathfrak{D}$-balanced evaluation functor $(M,F)\mapsto F(M)$. Finally, the 3-universal property of the relative 2-Deligne tensor product ensures that the diagram below commutes up to $\mathfrak{C}$-$\mathfrak{D}$-bimodule 2-natural equivalences: \begin{equation*}\begin{tikzcd}
\mathfrak{M} &  & \\
\mathfrak{M}\boxtimes_{\mathfrak{D}}\mathfrak{D} \arrow[d, "Id\boxtimes \eta"'] \arrow[u, "\simeq"] \arrow[rrd, "Id\boxtimes \{D\mapsto (-)\Box D\}", bend left=14]     &  & \\
{\mathfrak{M}\boxtimes_{\mathfrak{D}}\mathbf{Fun}_{\mathfrak{C}}(\mathfrak{M},\mathfrak{C})\boxtimes_{\mathfrak{C}}\mathfrak{M}} \arrow[d, "\epsilon\boxtimes Id"'] \arrow[rr, "\simeq"] &  & \mathfrak{M}\boxtimes_{\mathfrak{D}}\mathbf{End}_{\mathfrak{C}}(\mathfrak{M}) \arrow[lldd, "M\boxtimes F\mapsto F(M)", bend left=26] \\
\mathfrak{C}\boxtimes_{\mathfrak{C}}\mathfrak{M} \arrow[d, "\simeq"']                                                  &  & \\
\mathfrak{M}.  &  &                                        \end{tikzcd}\end{equation*} This prove one of the triangle identities. The other follows similarly, which shows that $\mathbf{Fun}_{\mathfrak{C}}(\mathfrak{M},\mathfrak{C})$ is indeed a left adjoint for $\mathfrak{M}$ in $\mathbf{CST2C}$. Likewise, one shows that $\mathbf{Fun}_{\mathfrak{D}}(\mathfrak{M},\mathfrak{D})$ is a right adjoint for the 1-morphism $\mathfrak{M}$. 

We now prove that any 2-morphism of $\mathbf{CST2C}$ whose source and target are locally separable compact semisimple bimodule 2-categories has a right and a left adjoint. Thanks to propositions 4.2.3 and 4.2.4 of \cite{D8}, it is enough to show that the underlying 2-functor has right and left adjoints, but this follows readily from corollary 3.2.3 of \cite{D5}. 

Finally, we show that every 3-morphism of $\mathbf{CST2C}$, that is bimodule 2-natural transfromation, which goes between two bimodule 2-functors whose source and target are separable and locally separable bimodule 2-categories has right and left adjoints. We begin by observing that every 2-natural transformation between locally separable compact semisimple 2-categories has right and left adjoints:\ This follows from theorem 3.4.3 of \cite{DSPS13}. Further, an adjoint 2-natural transformation to a bimodule 2-natural transformation inherits a canonical (op)lax bimodule structure compatible with the unit and counit of the adjunction. Thus, it only remains to prove that this (op)lax bimodule structure is actually strong. It is enough to do so for the left and right module structures separately. But, module 2-natural transformations between separable and locally separable module 2-categories have right and left adjoints as module 2-natural transformations thanks to theorem 5.1.2 of \cite{D8} and proposition 3.1.3 of \cite{D7}.
\end{proof}

\begin{Corollary}
Over a field of characteristic zero, every compact semisimple tensor 2-category is 4-dualizable.
\end{Corollary}

\begin{Remark}\label{rem:BrFus4Vect}
Let us assume that $\mathds{k}$ is an algebraically closed field of characteristic zero. The symmetric monoidal 4-category $\mathbf{BrFus}$ of braided multifusion 1-categories was introduced in \cite{BJS}, where it is shown that every object is fully dualizable. This 4-category is related to the symmetric monoidal 4-category $\mathbf{F2C}$ of multifusion 2-categories and finite semisimple bimodule 2-categories.\footnote{Following the higher condensation perspective of \cite{GJF}, we could also denote this symmetric monoidal 4-category by $\mathbf{4Vect}$.} Namely, in \cite{JF}, the construction of a symmetric monoidal 4-functor $\mathbf{Mod}:\mathbf{BrFus}\rightarrow \mathbf{F2C}$ is sketched. On objects, it is given by sending a braided multifusion 1-category $\mathcal{B}$ to the multifusion 2-category $\mathbf{Mod}(\mathcal{B})$. This sketch relies on the theory of higher condensations from \cite{GJF}, which has not yet been made fully rigorous. Nevertheless, as is explained in \cite{D9}, the 4-functor $\mathbf{Mod}$ would be an equivalence.\footnote{Up to inserting the obvious separability hypotheses, the characteristic zero hypothesis may be relaxed. We do not know whether this is also the case for the assumption that the field $\mathds{k}$ be algebraically closed.} Thus, theorem \ref{thm:fulldualizability} would follow from the main theorem of \cite{BJS}. Nevertheless, we emphasize that this approach does not bypass theorem \ref{thm:relative2Deligne}; It is still needed to argue that the symmetric monoidal 4-category $\mathbf{F2C}$ even exists.
\end{Remark}

\begin{Remark}
Dropping the separability and semisimplicity assumptions in theorem 3.2.4 of \cite{DSPS13}, one can still obtain partial dualizability results \cite[Thm. 3.2.2 \& Thm. 3.3.1]{DSPS13}. We expect that similar variants of our theorem \ref{thm:fulldualizability} hold, i.e.\ we expect that ``finite tensor 2-categories are 3-dualizable''. Under suitable connectedness hypotheses, this is essentially theorem 5.16 of \cite{BJS}. In particular, it would be interesting to understand whether ``every finite tensor 2-categories is Morita equivalent to a connected one''.
\end{Remark}

\begin{Remark}
In \cite{DSPS13}, a proof of the Radford isomorphism for finite tensor 1-categories was given using the fact that such categories have enough dualizability as objects of $\mathbf{TC}$, the Morita 3-category of finite tensor 1-categories. We expect that a similar analysis may be used to derive the categorified version of this result. More precisely, we expect that for any separable compact semisimple tensor 2-category, there is a monoidal 2-natural equivalence $(-)^{\sharp\sharp\sharp\sharp}\simeq Id$. We will not pursue this line of investigation further here. Nevertheless, for use below, we do collect the following related observation.
\end{Remark}

\begin{Lemma}\label{lem:functordualsagree}
Let $\mathfrak{C}$ and $\mathfrak{D}$ be locally separable compact semisimple tensor 2-categories, and let $\mathfrak{M}$ be a separable $\mathfrak{C}$-$\mathfrak{D}$-bimodule 2-category. Then, there is an equivalence of $\mathfrak{D}$-$\mathfrak{C}$-bimodule 2-categories \begin{equation*}\mathbf{Fun}_{\mathfrak{C}}(\mathfrak{M},\mathfrak{C})\simeq \mathbf{Fun}_{\mathfrak{D}}(\mathfrak{M},\mathfrak{D}).\end{equation*}
\end{Lemma}
\begin{proof}
This follows from lemma 1.4.4 of \cite{DSPS13} (see also remark 2.4.22 of \cite{L}) and the theorem above.
\end{proof}

\subsection{Invertibility}

Let $\mathds{k}$ be a perfect field. In the previous section, we have studied which objects and morphisms of the symmetric monoidal 4-category $\mathbf{CSST2C}$ are fully dualizable. We now examine which objects and 1-morphisms of this monoidal 4-category are invertible. We begin by carefully examining which 1-morphisms are invertible, i.e.\ which separable bimodule 2-categories witness a Morita equivalence.

At the decategorified level, i.e.\ for finite semisimple tensor 1-categories, invertibility of a finite semisimple bimodule 1-category admits various equivalent characterizations as is explained in proposition 4.2 of \cite{ENO2}. We prove the categorified version of this result below, which will be applied to study group-graded extensions of compact semisimple tensor 2-categories in \cite{D11}. Before stating this proposition, we need to recall the following definition.

Let us fix $\mathfrak{C}$ and $\mathfrak{D}$ two locally separable compact semisimple tensor 2-categories. A compact semisimple left $\mathfrak{C}$-module 2-category $\mathfrak{M}$ is \textit{faithful} provided that every simple summand of the monoidal unit of $\mathfrak{C}$ has non-zero action on $\mathfrak{M}$. We say that a compact semisimple $\mathfrak{C}$-$\mathfrak{D}$-bimodule 2-category is faithful if it is faithful both as a left and as a right module 2-category.

\begin{Proposition}\label{prop:invertibilitycharacterization}
Let $\mathfrak{M}$ be a faithful locally separable compact semisimple $\mathfrak{C}$-$\mathfrak{D}$-bimodule 2-category. The following are equivalent:
\begin{enumerate}
    \item The canonical monoidal 2-functor $\mathfrak{D}^{mop}\rightarrow \mathbf{End}_{\mathfrak{C}}(\mathfrak{M})$ is an equivalence.
    \item The canonical monoidal 2-functor $\mathfrak{C}\rightarrow\mathbf{End}_{\mathfrak{D}}(\mathfrak{M})$ is an equivalence.
    \item The 2-category $\mathfrak{M}$ is separable as a right $\mathfrak{D}$-module 2-category, and there is an equivalence of $\mathfrak{C}$-$\mathfrak{C}$-bimodule 2-categories \begin{equation*}\mathfrak{M}\boxtimes_{\mathfrak{D}}\mathbf{Fun}_{\mathfrak{D}}(\mathfrak{M},\mathfrak{D})\simeq \mathfrak{C}.\end{equation*}
    \item The 2-category $\mathfrak{M}$ is separable as a left $\mathfrak{C}$-module 2-category, and there is an equivalence of $\mathfrak{D}$-$\mathfrak{D}$-bimodule 2-categories \begin{equation*}\mathbf{Fun}_{\mathfrak{C}}(\mathfrak{M},\mathfrak{C})\boxtimes_{\mathfrak{C}}\mathfrak{M}\simeq \mathfrak{D}.\end{equation*}
    \item The separable $\mathfrak{C}$-$\mathfrak{D}$-bimodule 2-category $\mathfrak{M}$ defines an invertible 1-morphism from $\mathfrak{C}$ to $\mathfrak{D}$ in $\mathbf{CST2C}$.
\end{enumerate}
\end{Proposition}

\begin{Definition}
A faithful locally separable compact semisimple $\mathfrak{C}$-$\mathfrak{D}$-bimodule 2-category $\mathfrak{M}$ is invertible if it satisfies any, hence all, of the equivalent conditions of proposition \ref{prop:invertibilitycharacterization}.
\end{Definition}

\begin{proof}
That $1$ and $2$ are equivalent follows from corollary 5.4.4 of \cite{D8}. In particular, if either $1$ or $2$ is satisfied, then $\mathfrak{M}$ is separable as a $\mathfrak{C}$-$\mathfrak{D}$-bimodule 2-category.

Let us assume that $\mathfrak{M}$ is separable as a left $\mathfrak{C}$-module 2-category. It then follows from proposition \ref{prop:modulefunctorequivalence} that there is an equivalence \begin{equation*}\mathbf{Fun}_{\mathfrak{C}}(\mathfrak{M},\mathfrak{C})\boxtimes_{\mathfrak{C}}\mathfrak{M}\simeq \mathbf{End}_{\mathfrak{C}}(\mathfrak{M})\end{equation*} of $\mathfrak{D}$-$\mathfrak{D}$-bimodule 2-categories. In particular, $1$ implies $4$. Conversely, if we assume that $4$ holds, then we can consider the $\mathfrak{D}$-$\mathfrak{D}$-bimodule 2-functor \begin{equation*}\mathfrak{D}\rightarrow \mathbf{End}_{\mathfrak{C}}(\mathfrak{M})\simeq\mathfrak{D},\end{equation*} which is the composite of the canonical monoidal 2-functor with a pseudo-inverse of the 2-functor supplied by $4$. We claim that this composite is a $\mathfrak{D}$-$\mathfrak{D}$-bimodule 2-functor $\mathbf{E}:\mathfrak{D}\rightarrow\mathfrak{D}$ such that the image of the monoidal unit $I$ of $\mathfrak{D}$ is an invertible object of $\mathfrak{D}$ (as a monoidal 2-category). But, any such 2-functor $\mathbf{E}$ is automatically an equivalence, so that $1$ holds. In order to prove the claim, consider $\mathbf{F}:\mathfrak{D}\rightarrow \mathbf{End}_{\mathfrak{C}}(\mathfrak{M})$ an equivalence of $\mathfrak{D}$-$\mathfrak{D}$-bimodule 2-categories, and $J$ the object of $\mathfrak{D}$ such that $\mathbf{F}(J)=Id$. We have \begin{equation*}\mathbf{F}(I)\circ (-\Box J)\simeq \mathbf{F}(J)=Id\simeq (-\Box J)\circ \mathbf{F}(I)\end{equation*} in $\mathbf{End}_{\mathfrak{C}}(\mathfrak{M})$. This shows that $\mathbf{F}(I)$ is invertible as an object of $\mathbf{End}_{\mathfrak{C}}(\mathfrak{M})$ with inverse $(-\Box J)$. As duals are unique, we immediately find that $\mathbf{F}(I)\simeq (-\Box J^{\sharp})$, from which it follows that $J$ is an invertible object of $\mathfrak{D}$ as asserted.

Now, we assume that $\mathfrak{M}$ is separable as a right $\mathfrak{D}$-module 2-category. It follows from a variant of lemma \ref{lem:functordualsagree} that there is an equivalence of $\mathfrak{C}$-$\mathfrak{D}$-bimodule 2-categories \begin{equation*}\mathbf{Fun}_{\mathfrak{D}}(\mathbf{Fun}_{\mathfrak{D}}(\mathfrak{M},\mathfrak{D}),\mathfrak{D})\simeq \mathfrak{M}.\end{equation*} In particular, there is an equivalence of $\mathfrak{C}$-$\mathfrak{C}$-bimodule 2-categories \begin{equation*}\mathfrak{M}\boxtimes_{\mathfrak{D}}\mathbf{Fun}_{\mathfrak{D}}(\mathfrak{M},\mathfrak{D})\simeq \mathbf{End}_{\mathfrak{D}}(\mathbf{Fun}_{\mathfrak{D}}(\mathfrak{M},\mathfrak{D}))\simeq \mathbf{End}_{\mathfrak{D}}(\mathfrak{M}^*).\end{equation*} One can show that the assignment $F\mapsto \{M^*\mapsto \big(F^{*}(M)\big)^*\}$ gives a monoidal equivalence of 2-categories $\mathbf{End}_{\mathfrak{D}}(\mathfrak{M})\simeq \mathbf{End}_{\mathfrak{D}}(\mathfrak{M}^*)$. However, it is not quite compatible with the $\mathfrak{C}$-$\mathfrak{C}$-bimodule structures. Instead, the left $\mathfrak{C}$-module structure has to be twisted by the monoidal equivalence $(-)^{\sharp\sharp}:\mathfrak{C}\rightarrow \mathfrak{C}$. Nevertheless, this is enough to show that $2$ and $3$ are equivalent via an argument similar to that used to establish that $1$ and $4$ are equivalent.

Finally, by lemma \ref{lem:functordualsagree}, there is an equivalence of $\mathfrak{D}$-$\mathfrak{C}$-bimodule 2-categories \begin{equation*}\mathbf{Fun}_{\mathfrak{D}}(\mathfrak{M},\mathfrak{D})\simeq \mathbf{Fun}_{\mathfrak{C}}(\mathfrak{M},\mathfrak{C}),\end{equation*} so that $5$ is equivalent to $3$ and $4$ combined. This finishes the proof.
\end{proof}

Note that it follows from the definition that an invertible $\mathfrak{C}$-$\mathfrak{D}$-bimodule 2-category $\mathfrak{M}$ witnesses a Morita equivalence between $\mathfrak{C}$ and $\mathfrak{D}$ as defined in theorem 5.4.3 of \cite{D8}. Let us also point out that, over an algebraically closed field of characteristic zero, invertible bimodule 2-categories (in the sense of part (1) of proposition \ref{prop:invertibilitycharacterization}) have already appeared in definition 1.3.1 of \cite{DY23}. 

We now turn our attention to the invertible objects of $\mathbf{CSST2C}$, we begin by examining an example.

\begin{Example}
Let $\mathds{k}\subseteq\mathbb{K}$ be a finite Galois extension, with Galois group $G$. Then, we can consider $\mathbf{2Vect}_{\mathbb{K}}$ as a compact semisimple $\mathbf{2Vect}_{\mathds{k}}$-module 2-category. It is separable as the field extension is separable, so that it follows from \cite{D8} that the corresponding dual tensor 2-category $\mathfrak{K}$ is compact semisimple. One checks that the compact semisimple tensor 2-category $\mathfrak{K}$ is equivalent to $\mathbf{2Vect}_{\mathbb{K}}^{\sigma}(G)$, the compact semisimple tensor 2-category of $G$-graded 2-vector spaces over $\mathbb{K}$ with conjugation action of $G$ on $\Omega^2\mathbf{2Vect}_{\mathbb{K}}^{\sigma}(G)\cong\mathbb{K}$ given by the canonical action by field automorphisms. In particular, $\mathfrak{K}$ is a non-trivial compact semisimple tensor 2-category whose Morita equivalence class is trivial, and therefore also invertible.

More generally, if we only assume that $\mathds{k}\subseteq\mathbb{K}$ is a finite separable extension, we can still consider the Morita dual compact semisimple tensor 2-category $\mathfrak{K}$ of $\mathbf{2Vect}_{\mathds{k}}$ with respect to $\mathbf{2Vect}_{\mathbb{K}}$. However, it is somewhat trickier to describe in general. For instance, if we consider the extension $\mathbb{Q}\subseteq \mathbb{Q}(\sqrt[3]{2})$, then the underlying compact semisimple 2-category is \begin{equation*}\mathfrak{K}\simeq \mathbf{2Vect}_{\mathbb{Q}(\sqrt[3]{2})}\boxplus \mathbf{2Vect}_{\mathbb{Q}(\sqrt[3]{2}, e^{2\pi i/3})}.\end{equation*} This slightly generalized version of the construction also produces a non-trivial compact semisimple tensor 2-category whose Morita equivalence class is trivial, and therefore also invertible.
\end{Example}

\begin{Definition}
A separable compact semisimple tensor 2-category $\mathfrak{C}$ is (Morita) invertible if it is invertible as an object of $\mathbf{CST2C}$, or, more concretely, if $\mathfrak{C}\boxtimes\mathfrak{C}^{mop}$ is Morita equivalent to $\mathbf{2Vect}_{\mathds{k}}$.
\end{Definition}

As expected, invertibility can be checked at the level of the Drinfeld center (see \cite{Cr} for the original definition of the Drinfeld center of a monoidal 2-category, and \cite{D9} for a treatment in the context of fusion 2-categories over an algebraically closed field of characteristic zero).

\begin{Lemma}
Let $\mathfrak{C}$ be a separable compact semisimple tensor 2-category. Then, $\mathfrak{C}$ is invertible if and only if there is a monoidal equivalence $\mathscr{Z}(\mathfrak{C})\simeq\mathbf{2Vect}_{\mathds{k}}$.
\end{Lemma}
\begin{proof}
For any separable compact semisimple tensor 2-category $\mathfrak{C}$, there are equivalences of 2-categories \begin{equation*}\mathscr{Z}(\mathfrak{C})\simeq \mathbf{End}_{\mathfrak{C}\mathrm{-}\mathfrak{C}}(\mathfrak{C})\simeq\mathbf{Fun}_{\mathfrak{C}\boxtimes\mathfrak{C}^{mop}}(\mathfrak{C},\mathfrak{C}\boxtimes\mathfrak{C}^{mop})\boxtimes_{\mathfrak{C}\boxtimes\mathfrak{C}^{mop}}\mathfrak{C}.\end{equation*} The first equivalence is lemma 2.2.1 of \cite{D9}, and the second is proposition \ref{prop:modulefunctorequivalence}. In particular, the Drinfeld center as a compact-semisimple 2-category is invariant under Morita equivalence. Moreover, for any two separable compact semisimple tensor 2-categories $\mathfrak{C}$ and $\mathfrak{D}$, we have \begin{equation*}\mathscr{Z}(\mathfrak{C}\boxtimes\mathfrak{D})\simeq \mathscr{Z}(\mathfrak{C})\boxtimes\mathscr{Z}(\mathfrak{D})\end{equation*} as 2-categories. The forward implication now follows readily. Conversely, if $\mathscr{Z}(\mathfrak{C})\simeq\mathbf{2Vect}_{\mathds{k}}$, we have $\mathscr{Z}(\mathfrak{C}\boxtimes\mathfrak{C}^{mop})\simeq \mathbf{2Vect}_{\mathds{k}}$, so that $\mathfrak{C}\boxtimes\mathfrak{C}^{mop}$ is indecomposable, i.e.\ it cannot be written as the direct sum of two non-zero compact semisimple tensor 2-categories. This guarantees that the canonical left action of $\mathfrak{C}\boxtimes\mathfrak{C}^{mop}$ on $\mathfrak{C}$ is faithful, so that the left $\mathfrak{C}\boxtimes\mathfrak{C}^{mop}$-module 2-category $\mathfrak{C}$ witnesses a Morita equivalence between $\mathfrak{C}\boxtimes\mathfrak{C}^{mop}$ and $\mathscr{Z}(\mathfrak{C})$, and $\mathfrak{C}$ is indeed invertible. 
\end{proof}

The next lemma provides examples of invertible compact semisimple tensor 2-categories that are not necessarily Morita trivial. The related question of characterizing the invertible objects of $\mathbf{BrFus}$ was studied in \cite{BJSS} (see in particular theorem 4.2 therein for the complete answer over an algebraically closed field of characteristic zero).

\begin{Lemma}
Let $\mathcal{B}$ be a braided separable finite semisimple tensor 1-category. Then, the locally separable compact semisimple tensor 2-category $\mathbf{Mod}(\mathcal{B})$ is invertible if and only if $\mathcal{Z}_{(2)}(\mathcal{B})$, the symmetric center of $\mathcal{B}$, is $\mathbf{Vect}_{\mathds{k}}$.
\end{Lemma}
\begin{proof}
Assuming that $\mathcal{Z}_{(2)}(\mathcal{B})\simeq\mathbf{Vect}_{\mathds{k}}$, we find that the canonical braided functor $\mathcal{B}\boxtimes\mathcal{B}^{rev}\rightarrow\mathcal{Z}(\mathcal{B})$ is an equivalence. For instance, the proof given in proposition 8.6.3 of \cite{EGNO} does apply thanks to our hypotheses. Then, the proof of the lemma proceeds as in example 5.3.7 of \cite{D8}. Conversely, it follows from a slight generalization of theorem 4.10 of \cite{DN} that $\Omega\mathscr{Z}(\mathbf{Mod}(\mathcal{B}))\simeq \mathcal{Z}_{(2)}(\mathcal{B})$. In particular, if $\mathcal{Z}_{(2)}(\mathcal{B})$ is strictly larger than $\mathbf{Vect}_{\mathds{k}}$, then $\mathbf{Mod}(\mathcal{B})$ can not be invertible.
\end{proof}

\begin{Remark}
Over an algebraically closed field of characteristic zero, the converse of the above result also holds by lemma 4.1.3 of \cite{D9}. More precisely, every invertible fusion 2-category is Morita equivalent to $\mathbf{Mod}(\mathcal{B})$ with $\mathcal{B}$ a braided fusion 1-category with $\mathcal{Z}_{(2)}(\mathcal{B})\simeq\mathbf{Vect}_{\mathds{k}}$. Up to adding the necessary separability hypotheses, the proof also works over an arbitrary algebraically closed field. Yet more generally, over an arbitrary perfect field, one can show that every invertible compact semisimple tensor 2-category $\mathfrak{C}$ with $\Omega^2\mathfrak{C}\cong\mathds{k}$ is Morita equivalent to $\mathbf{Mod}(\mathcal{B})$ with $\mathcal{B}$ a braided separable tensor 1-category with $\mathcal{Z}_{(2)}(\mathcal{B})\simeq\mathbf{Vect}_{\mathds{k}}$. We do not know whether this fact still holds true without the assumption that $\Omega^2\mathfrak{C}\cong\mathds{k}$.\footnote{Nevertheless, let us point out that if a compact semisimple tensor 2-category $\mathfrak{C}$ over $\mathds{k}$ is actually linear over a proper field extension $\mathds{K}$ of $\mathds{k}$, then it can not be invertible as $\mathscr{Z}(\mathfrak{C})$ must then be linear over $\mathds{K}$.}
\end{Remark}

\subsection{Field Theoretic Interpretation}\label{sub:TFTs}

We discuss some applications of our results and related conjectures in the context of fully extended field theories in the spirit of \cite{DSPS13}. For simplicity, we take $\mathds{k}$ to be an algebraically closed field of characteristic zero. Recall that we write $\mathbf{F2C}$ for the symmetric monoidal 4-category of multifusion 2-categories, and finite semisimple bimodule 2-categories.

We write $\mathsf{Bord}^{fr}_4$ for the symmetric monoidal 4-category of 4-framed bordisms as constructed in \cite{CS}. Informally, this is the 4-category with objects 4-framed 0-dimensional manifolds, ..., 4-morphisms are (diffeomorphisms classes of) 4-framed manifolds with corners. We refer the reader to \cite[section 1.1]{DSPS13} for a review of the notion of $n$-framed manifolds. Now, for any symmetric monoidal 4-category $\mathscr{C}$, a 4-dimensional 4-framed fully extended topological field theory with value in $\mathscr{C}$ is a symmetric monoidal functor $\mathcal{F}:\mathsf{Bord}^{fr}_4\rightarrow\mathscr{C}$. By the cobordism hypothesis \cite{BD,L}, such a functor is completely characterised by its value on the positively 4-framed point. Further, the corresponding object of $\mathscr{C}$ is fully dualizable. More precisely, evaluation on the point induces a homotopy equivalence between the space of 4-dimensional 4-framed fully extended topological field theory with value in $\mathscr{C}$ and the spaces $(\mathscr{C}^{f.d.})^{\simeq}$ of fully dualizable objects in $\mathscr{C}$.

\begin{Corollary}
For every multifusion 2-category $\mathfrak{C}$, there is an essentially unique 4-dimensional 4-framed fully extended topological field theory $\mathcal{F}_{\mathfrak{C}}:\mathsf{Bord}^{fr}_4\rightarrow\mathbf{F2C}$, whose value on the positively 4-framed point is $\mathfrak{C}$.
\end{Corollary}

\noindent Naturally, the uniqueness statement in the above corollary refers back to the structure of the 4-category $\mathbf{F2C}$. More precisely, two multifusion 2-categories are equivalent as objects of $\mathbf{F2C}$ if and only if they are Morita equivalent in the sense of \cite{D8} as was shown in proposition \ref{prop:invertibilitycharacterization} above. But, it was established in theorem 4.2.2 of \cite{D9} that every fusion 2-category is Morita equivalent to a connected fusion 2-category. It is therefore overwhelmingly likely that all the field theories of the above corollary were already constructed in \cite{BJS}. However, as was discussed in remark \ref{rem:BrFus4Vect} above, making this precise would involve constructing a suitably functor $\mathbf{BrFus}\rightarrow \mathbf{F2C}$. Further, as we will see in corollary \ref{cor:fieldtheoreticdecomposition} below, the 4-category $\mathbf{F2C}$ is in some sense easier to study than $\mathbf{BrFus}$.

Let us say that a field theory is \textit{indecomposable} if it cannot be written as the direct sum of two non-zero field theories. Clearly, if $\mathfrak{C}$ is a fusion 2-category, the field theory $\mathcal{F}_{\mathfrak{C}}$ is indecomposable. More generally, $\mathcal{F}_{\mathfrak{C}}$ is indecomposable if and only if the multifusion 2-category $\mathfrak{C}$ is indecomposable, that is it cannot be written as the direct sum of two non-zero multifusion 2-categories.

It was proven in theorem 4.1.6 of \cite{D9} that every fusion 2-category is Morita equivalent to the 2-Deligne tensor product of a strongly fusion 2-category and an invertible fusion 2-category. On one hand, invertible fusion 2-categories are precisely the invertible objects of $\mathbf{F2C}$, that is the corresponding field theories are \textit{invertible}. It is known that every invertible fusion 2-category is of the form $\mathbf{Mod}(\mathcal{B})$ with $\mathcal{B}$ a non-degenerate braided fusion 1-category by lemma 4.1.3 of \cite{D9}. Further, it was shown in examples 5.3.8 and 5.4.6 of \cite{D8} that the notion of Morita equivalence between such fusion 2-categories corresponds precisely to the notion of Witt equivalence between non-degenerate braided fusion 1-categories introduced in \cite{DMNO}. The corresponding equivalence classes are classified by the so-called Witt group $\mathcal{W}itt$ of non-degenerate braided fusion 1-categories, whose structure is completely known \cite{DMNO, DNO, NRWZ, JFR}. On the other hand, the classification of bosonic strongly fusion 2-category was given in \cite{JFY}:\ They are all of the form $\mathbf{2Vect}^{\pi}(G)$ for a finite group $G$ a 4-cocycle $\pi$ for $G$ with coefficients in $\mathds{k}^{\times}$. The corresponding topological field theories are (framed versions of) Dijkgraaf-Witten theories. Finally, it will be shown in \cite{D11} that fermionic strongly fusion 2-categories are classified by a finite group $G$, together with homotopy-theoretic data. We refer to the topological field theories associated to strongly fusion 2-categories as \textit{generalized Dijkgraaf-Witten theories}.

\begin{Corollary}\label{cor:fieldtheoreticdecomposition}
Every indecomposable 4-dimensional 4-framed fully extended topological field theory $\mathcal{F}:\mathsf{Bord}^{fr}_4\rightarrow\mathbf{F2C}$ is the product of an invertible topologcial field theory and a generalized Dijkgraaf-Witten theory.
\end{Corollary}

\noindent Statements of this form originate in the literature on topological orders \cite{LKW, LW}. Furthermore, a version of the last result above was obtained in \cite[section V.D]{JF2} using the theory of higher condensations \cite{GJF}. We warn the reader that the above description is in general not an exact classification (see remark 4.1.7 of \cite{D9}, but also the upcoming work \cite{JF3}). However, it does become exact when it is restricted to bosonic fusion 2-categories by proposition 5.3.9 of \cite{D9}.

In many applications, it is desirable to consider bordisms equipped with other tangential structures, such as an orientation. We may consider the corresponding symmetric monoidal 4-category $\mathsf{Bord}^{or}_4$ of oriented bordisms as constructed in \cite{CS}. Then, a 4-dimensional oriented fully extended topological field theory is a symmetric monoidal functor $\mathsf{Bord}^{or}_4\rightarrow \mathscr{C}$. By the cobordism hypothesis \cite{L}, such field theories are classified by homotopy fixed points for the canonical action of $SO(4)$ on the space $(\mathscr{C}^{f.d.})^{\simeq}$. In the 3-dimensional case, this is expected to be related to the classical notion of sphericality for fusion 1-categories \cite[section 3.5]{DSPS13}. More precisely, it is conjectured that spherical fusion 1-categories carry a homotopy $SO(3)$ fixed point structure. In the 4-dimensional case, a notion of sphericality for fusion 2-categories was introduced in \cite{DR}, and, as was already alluded to in the previously cited work, it seems natural to make the following conjecture.

\begin{Conjecture}
Every spherical fusion 2-category admits an $SO(4)$ homotopy fixed point structure, and therefore provides an oriented fully extended field theory.
\end{Conjecture}

\noindent In fact, it seems likely that the notion of sphericality used in \cite{DR} ought to be weakened, and the above conjecture should be interpreted with this weaker structure in mind. Moreover, we expect that the state-sum construction for oriented 4-manifolds using spherical fusion 2-categories developed in \cite{DR} agrees with the value of the corresponding oriented fully extended topological field theory on 4-dimensional manifolds. Finally, we expect that every oriented fully extended topological field theory $\mathsf{Bord}^{or}_4\rightarrow\mathbf{F2C}$ is a \textit{Crane-Yetter theory}, i.e.\ arises from a braided fusion 1-category equipped with a ribbon structure (see \cite[section 1.5.1]{BJS} for a related discussion). In addition, corollary \ref{cor:fieldtheoreticdecomposition} should admit an appropriate oriented version.

\bibliography{bibliography.bib}

\end{document}